\documentclass[11pt]{amsart}
\usepackage{amsmath}

\theoremstyle{plain}
\newtheorem{thm}{Theorem}[section]
\newtheorem{lem}[thm]{Lemma}
\newtheorem{cor}[thm]{Corollary}
\newtheorem{prop}[thm]{Proposition}
\newtheorem{notitle}[thm]{ }

\newtheorem{ther}{Theorem}

\theoremstyle{definition}

\def \R {\mathbf{R}}
\def \Z {\mathbf{Z}}
\def \Sig{\Sigma}

\def \a {\alpha}
\def \b {\beta}
\def \g {\gamma}
\def \d {\delta}
\def \k {\kappa}
\def \lam {\lambda}

\def \l {\ell}
\def \o {\omega}
\def \s {\sigma}
\def \t {\tau}

\def \bd {\partial}
\def \x {\times}

\def \sw {\text{SW}}
\def \al {\text{A}}
\def \DD {\Delta}
\def \DN {\nabla}
\def \gr {\text{Gr}}
\def \Th {\Theta}
\def \C  {\mathcal{C}}

\begin{document}

\baselineskip.525cm

\title{Knots, Links, and 4-Manifolds}
\author[Ronald Fintushel]{Ronald Fintushel}
\address{Department of Mathematics, Michigan State University \newline
\hspace*{.375in}East Lansing, Michigan 48824}
\email{\rm{ronfint@math.msu.edu}}
\thanks{The first author was partially supported NSF Grant DMS9401032 and
the second
author by NSF Grant DMS9626330}
\author[Ronald J. Stern]{Ronald J. Stern}
\address{Department of Mathematics, University of California \newline
\hspace*{.375in}Irvine,  California 92697}
\email{\rm{rstern@math.uci.edu}}
\maketitle

\section{Introduction\label{Intro}}

In this paper we  investigate the  relationship between
isotopy classes of knots and links in $S^3$ and the diffeomorphism types
of homeomorphic smooth $4$-manifolds. As a  corollary of this initial
investigation, we begin to uncover the surprisingly
rich structure of diffeomorphism types of manifolds homeomorphic to the K3
surface.

In order to state
our theorems we need to view the Seiberg-Witten invariant of a smooth
$4$-manifold as a multivariable (Laurent) polynomial. To do this, recall
that the
Seiberg-Witten
invariant  of a smooth closed oriented $4$-manifold
$X$ with $b_2 ^+(X)>1$ is an integer valued function which is defined on
the set of $spin ^{\, c}$ structures over $X$, ({\it{cf.}} \cite{W},
\cite{KM},\cite{Ko},\cite{T1}). In case $H_1(X,\Z)$ has no 2-torsion (which
will be the
situation in this paper) there is a natural identification of the $spin ^{\, c}$
structures of
$X$ with the characteristic elements of $H^2(X,\Z)$. In this case we view the
Seiberg-Witten invariant as
\[ \sw_X: \lbrace k\in H^2(X,\Z)|k\equiv w_2(TX)\pmod2)\rbrace
\rightarrow \Z. \]
The Seiberg-Witten invariant $\sw_X$ is a diffeomorphism invariant whose sign
depends on an
orientation of
$H^0(X,\R)\otimes\det H_+^2(X,\R)\otimes \det H^1(X,\R)$. If $\sw_X(\b)\neq
0$, then we call
$\b$ a {\it{basic class}} of $X$. It is a fundamental fact that the set of
basic classes is
finite. Furthermore, if $\b$ is a basic class, then so is $-\b$ with
\[\sw_X(-\b)=(-1)^{(\text{e}+\text{sign})(X)/4}\,\sw_X(\b)\] where
$\text{e}(X)$ is
the Euler number and $\text{sign}(X)$ is the signature of $X$.

Now let
$\{\pm \b_1,\dots,\pm \b_n\}$ be the set of nonzero basic classes for $X$.
For the purposes of
this paper we define the Seiberg-Witten invariant of $X$ to be the formal series
\[\sw_X = b_0+\sum_{j=1}^n b_j(\exp(\b_j)+(-1)^{(\text{e}+\text{sign})(X)/4}\,
\exp(-\b_j))\]
where $b_0=\sw_X(0)$ and $b_j=\sw_X(\b_j)$.  Letting
$t_j=\exp(\b_j)$, we have that the Seiberg-Witten invariant is the `symmetric'
Laurent polynomial
\[\sw_X = b_0+\sum_{j=1}^n
b_j(t_j+(-1)^{(\text{e}+\text{sign})(X)/4}\,t_j^{-1})\]
in the (formal) variables $t_1,\dots,t_n$.

We refer to any symmetric
Laurent polynomial $P(t)=a_0+\sum\limits_{j=1}^na_j(t^j+t^{-j})$ of one
variable with coefficient sum
$a_0+2\,\sum\limits_{j=1}^n a_j=\pm1$ as an
$A$-{\it{polynomial}}. If, in addition, $a_n=\pm 1$ we refer to $P(t)$ as a
{\it{monic}}
$A$-polynomial.

Let $X$ be any simply connected smooth
4-manifold with
$b^+>1$. We
define a {\it{cusp}} in $X$ to be a PL embedded 2-sphere of
self-intersection 0 with a single
nonlocally flat point whose neighborhood is the cone on the right-hand
trefoil knot. (This agrees
with the notion of a cusp fiber in an elliptic surface.) The regular
neighborhood  $N$ of a cusp in a
ª4-manifold is a {\it {cusp neighborhood}}; it is the manifold obtained by
performing 0-framed surgery
on a trefoil knot in the boundary of the 4-ball. Since the trefoil
knot is a fibered knot with a genus 1 fiber, $N$ is fibered by smooth tori with
one singular fiber, the cusp.
If $T$ is a smoothly embedded torus representing a {\em{nontrivial}}
homology class $[T]$, we say that $T$ is {\it{c-embedded}} if
$T$ is a smooth fiber in a cusp neighborhood $N$; equivalently, $T$ has two
vanishing cycles. Note that a c-embedded torus has self-intersection $0$.
We can now state our first theorem.

\begin{thm}\label{introthm} Let $X$ be any simply connected smooth
4-manifold with
$b^+>1$. Suppose that $X$
contains a
smoothly c-embedded torus $T$ with $\pi_1(X\setminus T)=1$. Then for any
$A$-polynomial
$P(t)$, there is a smooth 4-manifold $X_{P}$ which is homeomorphic to $X$
and has
Seiberg-Witten invariant
\[ \sw_{X_{P}}=\sw_X\cdot P(t)\]
where $t=\exp(2[T])$.
\end{thm}

The basic classes of $X$ were defined above to be elements of $H^2(X)$. To
make sense of
the statement of the theorem, we need to replace $[T]$ by its Poincar\'e
dual. Throughout this paper we allow ourselves to pass freely between
$H^2(X)$ and
$H_2(X)$ without further comment.

As a corollary to the construction of $X_P$ we shall show:

\begin{cor}\label{symthm} Suppose further that $X$ is symplectic and that $T$ is
symplectically embedded. If $P(t)$ is a monic $A$-polynomial,
 then $X_P$ can be constructed as a symplectic manifold.
\end{cor}

Using the work of Taubes [T1--T5] concerning the
nature of the Seiberg-Witten invariants of symplectic manifolds, we shall
also deduce:

\begin{cor}\label{cor3} If $P(t)$ is not monic, then $X_P$ does not admit a
symplectic
structure. Furthermore, if $X$  contains a surface $\Sigma_g$ of genus $g$
disjoint from
$T$ with $0\ne[\Sigma_g] \in H_2(X;\Z)$ and with $[\Sigma_g]^2 < 2-2g$ if
$g>0$, or
$[\Sigma_g]^2 \le 0$ if $g=0$, then
$X_P$ with the opposite orientation  does not admit a symplectic structure.
\end{cor}
As a corollary, we have many interesting new homotopy
K3 surfaces ({\it{i.e.}} manifolds homeomorphic to the K3 surface). In
particular, since
$\sw_{K3}=1$ we have:

\begin{cor} \label{K3thm} Any $A$-polynomial $P(t)$ can occur as the
Seiberg-Witten invariant of an irreducible homotopy K3 surface.  If $P(t)$
is not monic, then the homotopy {\rm{K3}} surface
 does not admit a symplectic structure with either orientation. Furthermore,
any monic $A$-polynomial can occur as the Seiberg-Witten invariant of a
symplectic homotopy K3 surface.
\end{cor}

In fact, there are three disjoint c-embedded tori $T_1$, $T_2$, $T_3$ in the K3
surface representing
distinct homology classes $[T_j]$, $j=1,2,3$ ({\it{cf.}} \cite{GM}).
Theorem~\ref{introthm}
then implies that the product of any three $A$-polynomials $P_j(t_j)$,
$j=1,2,3$, can occur
as the Seiberg-Witten invariant of a homotopy K3 surface $K3_{P_1P_2P_3}$ with
$t_j=\exp(2[T_j])$. Furthermore, if all three of the polynomials are monic,
the resulting
homotopy K3 surface can be constructed as a symplectic manifold. If any one
of the
$A$-polynomials $P_j(t_j)$ is not monic, then the resulting homotopy K3
surface admits no
symplectic structure.

A common method for constructing exotic manifolds is to perform log
transforms on c-embedded tori.  One might ask whether
these new homotopy K3 surfaces $K3_{P_1P_2P_3}$
can be constructed in this fashion. However, it is shown in
\cite{rat} that if $X'$ is the result of performing log transforms of
multiplicities
$p_1,\dots,p_n$ on c-embedded tori in the K3 surface, then
$\sw_{X'}(1,\dots,1)=
\pm p_1\cdots p_n$. However,
$\sw_{K3_{P_1P_2P_3}}(1,1,1)=P_1(1)P_2(1)P_3(1)=\pm 1$; so
$K3_{P_1P_2P_3}$ cannot be
built using log transforms in this way. (Note that the `$1$' above is $\exp(0)$
corresponding to the zero class of $H_2$.)

The $A$-polynomials appearing in Theorem~\ref{introthm} are
familiar to knot theorists. It is known that any $A$-polynomial
occurs as the Alexander polynomial $\DD_K(t)$ of some knot $K$ in $S^3$.
Conversely, the Alexander polynomial of a knot $K$ is an $A$-polynomial.
Furthermore,  if
the $A$-polynomial is monic then the knot can be constructed as a fibered knot,
and if $K$ is fibered, then
$\DD_K(t)$ is a monic $A$-polynomial.
Indeed, it is a knot
$K$ in $S^3$ which we use to construct $X_P$.

Consider a knot $K$ in $S^3$, and let $m$ denote a meridional circle to
$K$. Let $M_K$ be the 3-manifold obtained by performing $0$-framed surgery
on $K$. Then $m$ can also be viewed as a circle in $M_K$. In $M_K\x S^1$ we
have the
smooth torus
$T_m=m\x S^1$ of self-intersection
$0$. Since a neighborhood of $m$ has a canonical framing in $M_K$, a
neighborhood of
the torus
$T_m$ in $M_K\x S^1$ has a canonical identification with $T_m\x D^2$. Let
$X_K$ denote the fiber sum
\[ X_K=X\#_{T=T_m}(M_K\x S^1)=[X\setminus (T\x D^2)] \cup [(M_K\x S^1)\setminus
(T_m\x D^2)]\]
where $T\x D^2$
is a  tubular neighborhood of the torus T in the manifold $X$ (with
$\pi_1(X)=\pi_1(X\setminus T)=0$).
The two pieces are glued together so as to
preserve the homology class $[{\text{pt}}\x \bd D^2]$. This latter
condition does not, in general,
completely determine the diffeomorphism type of $X_K$ ({\it{cf.}} \cite{Gompf}).
We take $X_K$ to be
any manifold constructed in this fashion.  Because
$M_K$ has the homology of
$S^2\x S^1$ with the class of $m$ generating $H_1$, the complement $(M_K\x
S^1) \setminus
(T\x D^2)$ has the homology of $T^2\x D^2$. Thus $X_K$ has the same
homology (and
intersection pairing) as $X$. Furthermore, the class of $m$ normally
generates
$\pi_1(M_K)$; so $\pi_1(M_K\x S^1)$ is normally generated by the image of
$\pi_1(T)$. Since
$\pi_1(X\setminus F)=1$, it follows from Van Kampen's Theorem that $X_K$
is simply
connected. Thus $X_K$ is homotopy equivalent to $X$. It is conceptually
helpful to note
that $X_K$ is obtained from
$X$ by removing a neighborhood of a torus and replacing it with the
complement of the knot $K$ in $S^3$
crossed with $S^1$.
Also, in order to define Seiberg-Witten invariants, the oriented 4-manifold $X$
must also be
equipped with an orientation of $H^2_+(X;\R)$. The manifold $X_K$ inherits
an orientation
as well as an orientation of $H^2_+(X_K;\R)$ from $X$.

Let $[T]$ be the class in $H_2(X_K;\Z)$ induced by the torus $T$ in $X$, and
let $t=\exp(2[T])$. Our first main theorem, from which
Theorem~\ref{introthm} is an
immediate corollary, is:

\begin{thm} \label{knotthm} If the torus $T$ is c-embedded, then
the Seiberg-Witten invariant of  $X_K$
is
\[ \sw_{X_K}=\sw_{X}\cdot\DD_K(t). \]
\end{thm}

Let $E(1)$ be the rational elliptic surface with elliptic fiber $F$. If $T$ is a
smoothly embedded self-intersection $0$ torus in $X$, then
$E(1)\#_{F=T}X_{K}=(E(1)\#_{F=T}X)\#_{T=T_m}(M_K\x S^1)$, and in
$E(1)\#_{F=T}X$, the torus $T=F$ is c-embedded. Thus we have a slightly more
general result:
\begin{cor} \label{strongknotthm} Let $X$ be any simply connected smooth
4-manifold with
$b^+>1$. Suppose that $X$
contains a
smoothly embedded torus $T$ of self-intersection $0$ with
$\pi_1(X\setminus T)=1$ and representing a nontrivial homology class $[T]$.
Then
\[ \sw_{E(1)\#_{F=T}X_{K}}=\sw_{E(1)\#_{F=T}X}\cdot \DD_K(t).\]
\end{cor}

More can be said  if $K$ is a fibered knot. Consider the normalized
Alexander polynomial $\al_K(t)=t^{d}\DD_K(t)$,
where $d$ is the degree of $\DD_K(t)$. If $K$ is a
fibered knot in $S^3$
with a punctured genus $g$ surface as fiber, then $M_K$ fibers over the
circle with a closed genus
$g$ surface $\Sigma_g$ as fiber. Thus
$M_K\x S^1$  fibers over $S^1\x S^1$ with $\Sigma_g$ as fiber and with $T_m=m\x
S^1$ as section.  It is a theorem of Thurston \cite{Th} that such a
4-manifold has
a symplectic structure with symplectic section $T_m$. Thus, if $X$ is a
symplectic $4-$manifold with a symplectically embedded torus with
self-intersection $0$, then  $X_K=X\#_{T=T_m}(M_K\x S^1)$ is symplectic since it
can be constructed as a symplectic fiber sum \cite{Gompf}. As a corollary
to Theorem~\ref{knotthm} and the theorems of Taubes relating the Seiberg-Witten
and Gromov invariants of a symplectic
$4-$manifold \cite{T3,T4} we have:

\begin{cor} \label{symcor} Let $X$ be a symplectic $4$-manifold with
$b^+>1$ containing
 a symplectic c-embedded torus $T$. If
$K$ is a fibered knot, then $X_K$ is a symplectic $4$-manifold whose Gromov
invariant is
\[\gr_{X_K}=\gr_X\cdot\al_K(\t)\]
where $\t=\exp([T])$.
\end{cor}
\begin{proof} The homology $H_2(M_K\x S^1)$ is generated by the classes of
the symplectic
curves $T_m$ and $\Sig_g$; so the canonical class of $M_K\x S^1$ has the form
$\k_{M_K\x S^1}=a[T_m]+b[\Sig_g]$. Applying the adjunction formula (and using
$[T_m]^2=[\Sig_g]^2=0$ and $[T_m]\cdot[\Sig_g]=1$) gives $b=0$ and
$a=2g-2$. But note
that the degree of $\DD_K(\t)=a_0+\sum\limits_{n=1}^da_n(\t^n+\t^{-n})$ is
$d=g$. Hence
$\k_{M_K\x S^1}=(2d-2)[T_m]$. This means
that the canonical class of the symplectic structure on $X_K$ is
$\k_{X_K}=\k_X+\k_{M_K\x S^1}+2[T]=\k_X+2d[T]$. Taubes' theorem now implies
that
for any $\a\in H_2(X_K)$, the coefficient of $\exp(\a)$ in $\gr_{X_K}$ is:
\begin{multline*}
\gr_{X_K}(\a)=\sw_{X_K}(2\a-\k_{X_K})=\sw_{X_K}(2\a-\k_X-2d[T])\\
    =\sum_{n=-d}^d a_n\,\sw_X(2\a-\k_X-2(d+n)[T])=\sum_{n=-d}^d
a_n\gr_X(\a-(d+n)[T]),
\end{multline*}
and this is the coefficient of $\exp(\a)$ in $\gr_X\cdot\al_K(\t)$.
\end{proof}

Of course, if $X$ is simply connected and $\pi_1(X\setminus T)=1$, then $X_K$ is
homeomorphic to $X$.  This implies Corollary~\ref{symthm}. As  corollary of the
initial work of
Taubes \cite{T1,T2} on the Seiberg-Witten invariants of symplectic
manifolds and the
adjunction inequality, we have the following corollary which also implies
Corollary~\ref{cor3}.

\begin{cor} \label{nonsymthm}  If $\DD_K(t)$ is not monic, then $X_K$ does
not admit a
symplectic structure. Furthermore, if
$X$  contains a surface $\Sigma_g$ of genus $g$ disjoint from
$T$  with $0 \ne[\Sigma_g]
\in H_2(X;\Z)$ and with $[\Sigma_g]^2 < 2-2g$ if $g>0$ or $[\Sigma_g]^2 < 0$ if
$g=0$, then $X_K$ with the opposite orientation  does not admit a symplectic
structure.
\end{cor}
\begin{proof} Suppose that $X_K$ admits a symplectic structure with
symplectic form
$\o$ and canonical class $\k$. Taubes has
shown that $\sw_{X_K}(\k)=\pm1$ and that if $k$ is any other
basic class then $|k\cdot\o|< \k\cdot\o$.

Since $\sw_{X_K}=\sw_{X}\cdot\DD_K(t)$ and since
$\DD_K(t)=a_0+\sum\limits_{n=1}^da_n(t^n+t^{-n})$ is a polynomial in one
variable
$t=\exp(2[T])$, any nontrivial basic class, and in particular,
the canonical class, is of the form
$\k=\a+2n[T]$ where $\sw_X(\a)\ne0$ and $|n|\le d$. Let $m$ be the maximum
integer satisfying $\sw_X(\a+m[T])\ne0$. Note that  $m\ge0$. Set
$\b=\a+m[T]$.
Because of the maximality of
$m$, we have $\sw_{X_K}(\b+2d[T])=a_d\cdot\sw_X(\b)\ne0$.

Replacing $[T]$ with $-[T]$ allows us to assume that $[T]\cdot\o\ge0$.
First assume that $[T]\cdot\o>0$. Then
\[(\b+2d[T])\cdot\o=\k\cdot\o+(m+2(d-n))[T]\cdot\o\ge\k\cdot\o\]
because $m\ge0$ and $d-n\ge0$, and equality occurs only if $m=0$ and $n=d$.
But a
strict inequality contradicts Taubes' theorem, thus $\a=\b$ and $n=d$; so
$\k=\b+2d[T]$.
This means that $\pm1=\sw_{X_K}(\k)=a_d\cdot\sw_X(\b)$; so $a_d=\pm1$,
{\it{i.e.}} $\DD_K$
is monic.

If $[T]\cdot\o=0$,
\[(\b+2d[T])\cdot\o=(\a+2n[T])\cdot\o=\k\cdot\o,\]
which means that $\k=\b+2d[T]$, and again we see that $a_d=\pm1$.

Finally, if any manifold $Y$  contains a homologically nontrivial surface
$\Sigma_g$ of
genus $g$ with $[\Sigma_g]^2 > 2g-2$, then, if $g>0$, it
follows from the adjunction inequality \cite{KM,MST} that the Seiberg-Witten
invariants of $Y$ vanish, and hence $Y$ does not admit a symplectic
structure \cite{T1}. If $g=0$ and $[\Sigma_g]^2 \ge 0$ one can also show
that the
Seiberg-Witten invariants of $Y$ must vanish \cite{FS,Ko1}.
\end{proof}

The authors are unaware of any simply connected smooth oriented
$4-$manifold $X$ with
$b^+>1$ and $SW_X\ne 0$ which does not contain
either a sphere with self-intersection $-2$ or a torus with self-intersection
$-1$.

The techniques used in proving Theorem~\ref{knotthm} generalize to the more
general setting of links. Let $L=\{K_1,\dots,K_n\}$ be an ($n\ge 2$)-component
ordered link in $S^3$ and suppose that $X_j$, $j=1,\dots,n$, are
simply connected smooth $4$-manifolds with $b^+\ge 1$ and each
containing a
smoothly embedded torus $T_j$ of self-intersection $0$ with
$\pi_1(X_j\setminus T_j)=1$ and representing a nontrivial homology class
$[T_j]$. If
$(\ell_j,m_j)$ denotes the standard longitude-meridian pair for the knot
$K_j$, let
\[\alpha_L:\pi_1(S^3\setminus L)\to \Z \]   denote the homomorphism
characterized by the
property that
$\alpha_L(m_j)=1$ for each
$j=1,\dots,n$. Now, mimicking the knot case above, let $M_L$ be the 3-manifold
obtained by performing
$\alpha_L(\ell_j)$ surgery on each component $K_j$ of $L$. (The surgery
curves form the
boundary of a Seifert surface for the link.) Then, in
$M_L\x S^1$ we have  smooth tori $T_{m_j}=m_j\x S^1$ of self-intersection $0$
and we can construct
the $n-$fold fiber sum
\[ X(X_1,\dots X_n;L)=(M_L\x S^1)\#_{T_j=T_{m_j}}\coprod\limits_{j=1}^n X_j\]
Here, the fiber sum is performed using the natural
 framings $T_{m_j} \x D^2$ of the neighborhoods of $T_{m_j}=m_j\x S^1$ in
$M_L\x S^1$
and the neighborhoods
$T_j
\x D^2$ in each $X_j$ and glued together so as to preserve the homology classes
$[{\text{pt}}\x
\bd D^2]$. As in the knot case, it follows from Van Kampen's Theorem that
$X(X_1,\dots X_n;L)$ is simply connected. Furthermore, the signature and Euler
characteristic ({\it{i.e.}} the rational homotopy type) of
$X(X_1,\dots X_n;L)$
depend only on the rational homotopy type of the $X_i$ and the number of
components in
the link $L$. In the special case that all the $X_j$ are the same manifold
$X$ we denote
the resulting construction by $X_L$. It is conceptually helpful to note
that $X(X_1,\dots X_n;L)$ is obtained
from
the disjoint union of the $X_j$ by removing a neighborhood of the tori $T_j$ and
replacing them with the complement of the link
$L$ in
$S^3$ crossed with
$S^1$.

If $\DD_L(t_1,\dots,t_n)$ denotes the symmetric multivariable Alexander
polynomial of the $n\ge 2$ component link $L$, and $E(1)$ denotes the
rational elliptic
surface, our second theorem is:

\begin{thm} \label{linkthm} The Seiberg-Witten invariant of
$E(1)_L$ is
\[
\sw_{E(1)_L}=\
\DD_L(t_1,\dots,t_n)
\] where $t_j=\exp(2[T_j])$.
\end{thm}

As a corollary we shall show:
\begin{cor} \label{linkcor} The Seiberg-Witten invariant of
$X(X_1,\dots X_n;L)$ is
\[
\sw_{X(X_1,\dots
X_n;L)}=\DD_L(t_1,\dots,t_n)\cdot\prod_{j=1}^n\sw_{E(1)\#_{F=T_j}X_j}
\] where $t_j=\exp(2[T_j])$.
\end{cor}
\noindent As before, the work of Taubes \cite{T3,T4} implies that if $L$ is a
fibered link, and each $(X_j,T_j)$ is a symplectic pair, then
$X(X_1,\dots X_n;L)$ is a symplectic 4-manifold whose Gromov invariant is
\[
\gr_{X(X_1,\dots
X_n;L)}=\al_K(t_1,\dots,t_n)\cdot\prod_{j=1}^n\gr(E(1)\#_{F=T_j}X_j).\]

Two  words of caution are in order here. First, there may be relations in
$X(X_1,\dots
X_n;L)$ among the homology classes $[T_j]$ represented by the tori
$T_j=T_{m_j}$. These
relations are determined by the linking matrix of the link $L$. In
particular, if all the
linking numbers are zero, then the $[T_j]$ are linearly independent. At the
other
extreme, if the $n$-component link is obtained from the Hopf link by
pushing off one
component $(n-2)$ times, then all the $[T_j]$ are equal. Second, the
ordering of the
components of the link can affect these relations.

If $L$ is a two component link with odd linking number, then $E(1)_L$ is a
homotopy
$K3$-surface and there are many interesting new polynomials which are not
products of
$A$-polynomials ({\it{cf.}}
\cite{Hill}) that can occur in this way as
the Seiberg-Witten invariants of a homotopy $K3$-surface.

The first examples of nonsymplectic simply connected irreducible smooth
$4$-manifolds
 were constructed by Z. Szabo \cite{S1}. These manifolds $X(k)$ ($k\in\Z,
k\ne 0,\pm1$) can
be shown to be diffeomorphic to $E(1)_{W(k)}$, where $W(k)$ is the $2$-component
$k$-twisted Whitehead link (see Figure 1)
with Alexander polynomial
$k(t_1^{1/2}-t_1^{-1/2})(t_2^{1/2}-t_2^{-1/2})$.
By
Theorem~\ref{linkthm} (and by the computation first given in \cite{S1})
\[\sw_{X(k)}=
k(t_1^{1/2}-t_1^{-1/2})(t_2^{1/2}-t_2^{-1/2}).\]
Thus by Taubes' Theorem \cite{T1}, $X(k)$ does not admit a symplectic structure
with either orientation (since $X(k)$ contains spheres with self-intersection
$-2$). Note also that for
$k=\pm 1$ the
$k-$twisted Whitehead link is fibered; so
$X(\pm1)$ is, in fact, symplectic.

The first examples of
nonsymplectic homotopy $K3$-surfaces were constructed by the authors. These
manifolds $Y(k)$ can be shown to be diffeomorphic to $K3_{T(k)}$
where $T(k)$ is the
$k$-twist knot (see Figure 1)
with Alexander
polynomial $kt-(2k+1)+kt^{-1}$. By Theorem~\ref{knotthm}
\[\sw_{Y(k)}=kt-(2k+1)+kt^{-1}\] and, by Corollary~\ref{nonsymthm}, if $k
\ne 0, \pm1$,
$Y(k)$
does not admit a symplectic structure with either orientation. Again, for
$k=\pm 1$ the
$k-$twist knot is fibered; so
$Y(\pm1)$ is symplectic.

\vspace{-.1in}
\centerline{\unitlength 1cm
\begin{picture}(9.25,4.5)
\put (.75,2.5){\oval(1.5,1)[l]}
\put (.75,2.5){\oval(2.5,2)[l]}
\put (1.75,2.5){\oval(1.5,1)[r]}
\put (1.75,2.5){\oval(2.5,2)[r]}
\put (1.25,3.6){\oval(1.25,.5)[t]}
\put (.75,3.25){\oval(.5,.5)[rt]}
\put (1.75,3.25){\oval(.5,.5)[lt]}
\put (1.25,3.4){\oval(1.25,.5)[b]}
\put (.75,3.075){\oval(.5,.15)[rb]}
\put (1.75,3.075){\oval(.5,.15)[lb]}
\put (.75,1.3){\framebox(1,.9)}
\put (.85,1.85){\Small{$k$ RH-}}
\put (.85,1.5){\Small{twists}}
\put (-1 ,.9){\small{$W(k)=k$-twisted Whitehead link}}
\put (6.75,2.5){\oval(1.5,1)[l]}
\put (6.75,2.5){\oval(2.5,2)[l]}
\put (7.75,2.5){\oval(1.5,1)[r]}
\put (7.75,2.5){\oval(2.5,2)[r]}
\put (6.75,3.49){\line(1,0){.25}}
\put (6.75,2.99){\line(1,0){.25}}
\put (7.5,3.49){\line(1,0){.25}}
\put (6.95,3.175){\oval(1,.65)[rt]}
\put (7.6,3.325){\oval(1,.65)[lb]}
\put (7.5,2.99){\line(1,0){.25}}
\put (6.75,1.3){\framebox(1,.9)}
\put (6.85,1.85){\Small{$k$ RH-}}
\put (6.85,1.5){\Small{twists}}
\put (5.9,.9) {\small{$T(k) =k$-twist knot}}
\put (3.9,.25){Figure 1}
\end{picture}}

\vspace{-.05in}
Our next task is to  prove Theorem~\ref{knotthm} and Theorem~\ref{linkthm}.
The proof of
Theorem~\ref{knotthm} is constructive and gives an algorithm which relates the
Seiberg-Witten invariants of $X_K$ with those of $X$ by performing a series
of topological
log transforms on nullhomologous tori in $X_K$ which reduce it to $X$. This
turns out
to be the same algorithm used to compute the Alexander polynomial
$\DD_K(t)$. This proof  relies upon   important analytical work of
 Morgan, Mrowka, and Szabo \cite{MMS} ({\it{cf.}} \cite{S1,S2}) and  Taubes
\cite{T5}
regarding the effect on the Seiberg-Witten invariants of removing
neighborhoods of tori and
sewing in manifolds with nonnegative scalar curvature.
 These we present in
Section~\ref{backgd} in the form of gluing theorems. The proof of
Theorem~\ref{linkthm}
can take one of two routes. The first is to utilize the algorithm provided
by Conway
\cite{C} to compute the Alexander polynomial (more precisely the potential
function)
of a link. The proof then proceeds in a (tedious) manner,  similar  to the
proof of
Theorem~\ref{knotthm}. However, we choose to present a more direct proof
by  showing that the Seiberg-Witten
invariants for the manifolds
$E(1)_L$ satisfy the axioms for the Alexander polynomial of a link as
provided by
Turaev \cite{Turaev}. The advantage of this proof, aside from its brevity,
is that it
isolates the required gluing properties of the Seiberg-Witten invariants
and perhaps lays
the foundation for determining the axioms for an appropriate gauge theory
which may expose
the other, more sophisticated, knot and link invariants.
In the final section we discuss examples with $b^+=1$ which are given by
our construction.

It was Meng and Taubes \cite{MT} who first discovered the relationship between
Seiberg-Witten type invariants and the Alexander polynomial. In \cite{MT} they
defined a 3-manifold invariant by dimension-reducing the Seiberg-Witten
invariants, and they showed that this 3-manifold invariant was related to the
Milnor torsion. We fell upon Theorems~\ref{knotthm} and
\ref{linkthm} by attempting to understand our above mentioned constructions of
nonsymplectic homotopy $K3$-surfaces.

We end this introduction with three items. First, we  conjecture that if $K$
and $K'$ are
two distinct knots  (or
$n$-component links) then $X_K$ is diffeomorphic to $X_{K'}$ if and only if
$K$ is isotopic to $K'$. Second, we wish to thank Jim Bryan, Bob Gompf,
Elly Ionel, Dieter
Kotschick, Wladek Lorek, Dusa McDuff, Terry Lawson, Tom Parker, and Cliff
Taubes for
useful conversations.
Finally, we wish to make it clear that the contributions of the present
paper are of a
purely topological nature. The gauge theoretic input to our theorems is due to
Morgan, Mrowka, and Szabo and to Taubes.

\section{Background for the proofs of Theorems \ref{knotthm} and \ref{linkthm}
\label{backgd}}

In this section we shall survey the recent gluing theorems of
Morgan, Mrowka, and Szabo \cite{MMS} ({\it{cf.}} \cite{S1,S2}) and Taubes
\cite{T5}
({\it{cf.}}
\cite{MT}) which are used in the proof of Theorem~\ref{knotthm}. Also we
shall review some
of the work of W. Brakes and J. Hoste on
`sewn-up link exteriors' which will be used in our constructions.

The context for the first of the gluing results is as follows. We are given
a smooth 4-manifold $X$ with $b_X^+>1$ and  with an embedded torus T which
represents a
nontrivial
homology class $[T]$ of self-intersection $0$ in $H_2(X;\Z)$. Any
Seiberg-Witten
basic class
$\a\in H_2(X;\Z)$ ({\it{i.e.}} any $\a$ with
$\sw_X(\a)\ne 0$) must be
orthogonal to the homology class $[T]$ since the adjunction
inequality states that
$0\ge [T]^2+|\a\cdot[T]|$. The {\it{relative Seiberg-Witten invariant}}
$\sw_{(X;T)}$ is formally
defined to be
\[ \sw_{(X;T)}=\sw_{X\#_{T=F}E(1)} \] where $E(1)$ is
the rational elliptic surface
with smooth elliptic fiber $F$.

It is an interesting consequence  of the gluing theorems of
 Morgan-Mrowka-Szabo \cite{MMS} and Taubes \cite{T5} that

\begin{thm}\label{eiglue} Suppose that $b_X^+>1$ and that the torus $T$
is c-embedded. Then
\[ \sw_{(X;T)}={\sw_X}\cdot(t^{1/2}-t^{-1/2}) \] where $t=\exp(2[T])$.
\end{thm}

Note that $E(1)$ has a metric of positive scalar curvature. A much more
general (and
difficult to prove) gluing theorem is:

\begin{thm} [Morgan, Mrowka, and Szabo \cite{MMS}] \label{fiberglue} In the
situation above
\[\sw_{X_1\#_{T_1=T_2}X_2}=\sw_{(X_1;T_1)}\cdot \sw_{(X_2;T_2)}.\]
\end{thm}

The fact that Corollary~\ref{linkcor} follows from Theorem~\ref{linkthm} is now
an easy consequence of the gluing
theorems Theorem~\ref{eiglue} and Theorem~\ref{fiberglue}; for note that
$X(X_1,\dots,X_n;L)=X(X_2,\dots,X_n;L)\#_{T_{m_1}=T_1}X_1$. Thus
Theorem~\ref{eiglue}
implies that
\[\sw_{X(X_1,\dots,X_n;L)}=\sw_{X(E(1),X_2,\dots,X_n;L)}\cdot\sw_{X_1\#_{T_1
=F}E(1)},\]
and continuing inductively completes the argument. This
is the only time we shall need to use the general gluing theorem
\eqref{fiberglue}.

For our proof of Theorem~\ref{knotthm} we will form an `internal fiber
sum'. For this
construction suppose that we have a pair of c-embedded tori $T_1$, $T_2$.

In our manifold $X$ with $b_X^+>1$
we formally define the
relative Seiberg-Witten invariant $\sw_{(X;T_1,T_2)}$  to be
\[ \sw_{(X;T_1,T_2)}={\sw_X}\cdot(\t_1-\t_1^{-1})\cdot(\t_2-\t_2^{-1}). \]
where
$\t_j=\exp([T_j])$. We construct the {\it{internal fiber sum}}
$X_{T_1,T_2}$ by identifying the boundaries of neighborhoods of the
$T_i$, again preserving  the homology classes $[{\text{pt}}\x\bd
D^2]$.

The first gluing theorem we need
for the proof of Theorem~\ref{knotthm} is:

\begin{thm} [Morgan, Mrowka, and Szabo \cite{MMS}, Taubes
\cite{T5}({\it{cf.}}\cite{MT})]
\label{gluing1} The internal fiber sum $X_{T_1,T_2}$ has Seiberg-Witten
invariant
\[ \sw_{X_{T_1,T_2}}=\sw_{(X;T_1,T_2)}|_{\t_1=\t_2}.\]
\end{thm}

The other gluing result we need concerns generalized log transforms on
nullhomologous tori. Let $p$ and
$q$ be relatively prime nonzero integers (or $(1,0)$ or $(0,1)$). If
$T$ is any embedded self-intersection $0$ torus in a 4-manifold $Y$ with tubular
neighborhood
\[ N=T\x D^2=S^1\x S^1\x D^2,\]
let $\varphi=\varphi_{p,q}$ be the diffeomorphism $S^1\x S^1\x\bd D^2\to\bd
N $ given by
\[ \varphi(x,y,z)=(x,y^sz^q,y^rz^p),  \hspace{.25in}
\det\begin{pmatrix} p&q\\r&s \end{pmatrix}=-1. \]
The manifold
\[ Y(p/q) = (Y\setminus N)\cup_{\varphi} (S^1\x S^1\x D^2) \]
is called the {\it {(generalized) $(p/q)$-log transform of $Y$ along $T$}}.
Our notation and
terminology are incomplete since the splitting $T=S^1\x S^1$ is necessary
information.
Throughout this paper, whenever a log transform is performed on a torus
$T$, there will be a natural identification $T=S^1\x S^1$ and we always perform
the transform with
$\varphi_{p,q}$ in the
coordinates $S^1\x S^1\x D^2$ as above.

We shall need to study the situation where a log transform is performed on
a nullhomologous
torus in  $Y$.  If $T$ is such a torus, then in $Y(0/1)$ there appears a new
2-dimensional
homology class which is represented by the torus
\[ T_0 = S^1\x S^1\x{\text{pt}}\subset  S^1\x S^1\x D^2\subset
(X\setminus N)\cup_{\varphi} (S^1\x S^1\x D^2)=Y(0/1).\]
(The old torus pushed to $\bd(Y\setminus N)$ is now
$\varphi(S^1\x{\text{pt}}\x\bd D^2)$.)
Notice that each homology class $\a\in H_2(Y;\Z)$ may be viewed as a homology
class in each $Y(p/q)$.

\begin{thm}[Morgan, Mrowka, and Szabo \cite{MMS}, Taubes \cite{T5} ({\it{cf.}}
\cite{MT})]\label{log} Let
$Y$ be a smooth 4-manifold with $b^+\ge3$, and suppose that $Y$ contains a
nullhomologous
torus $T$ with tubular neighborhood $N=T\x D^2=S^1\x S^1\x D^2$. Let $\t$
be the homology
class of $T_0$ in $Y(0/1)$. Then for each characteristic homology class
$\a\in H_2(Y;\Z)$,
\[\sw_{Y(p/q)}(\a)= p\sw_Y(\a) +
q\sum\limits_{i=-\infty}^{\infty}\sw_{Y(0/1)}(\a+2i\t).\]
\end{thm}

The sum in the above formula reduces to a single term in all the situations
which are encountered in the proofs of Theorem~\ref{knotthm} and
Theorem~\ref{linkthm}. Specifically,

\begin{notitle}\label{log+} In Theorem~\ref{log} suppose that there is a
torus in $Y(0/1)$ which is disjoint from $T$ and
which represents a homology class $\s$ of self-intersection $0$ whose
intersection number with $\t$ in $Y(0/1)$ is $1$. Suppose furthermore that
$\s\cdot\a=0$ for all $\a\in H_2(Y)\subset H_2(Y(0/1))$. Then
\[\sw_{Y(p/q)}= p\sw_Y + q\sw_{Y(0/1)}.\]
\end{notitle}
\begin{proof} Note that $H_2(Y(0/1))=H_2(Y)\oplus H(\s,\t)$ where $
H(\s,\t)$ is a
hyperbolic pair. If
$\a\in H_2(Y;\Z)$ satisfies
$\sw_{Y(0/1)}(\a+2i\t)\ne0$, the adjunction
inequality implies that
\[0\ge\s^2+|(\a+2i\t)\cdot\s|=|2i|.\]
Thus $i=0$, and
$\sw_{Y(p/q)}(\a)= p\sw_Y(\a) + q\sw_{Y(0/1)}(\a)$.
Since for $p\ne0$, each $\a\in H_2(Y(p/q))$ arises from a class in
$H_2(Y)$, the lemma
follows.
\end{proof}

We next wish to describe a method for constructing 3-manifolds which was first
studied by W. Brakes \cite{WB} and extended by J. Hoste \cite{Hoste}. Let
$L$ be a link in $S^3$ with two {\it{oriented}}\/ components $C_1$ and
$C_2$. Fix tubular neighborhoods $N_i\cong S^1 \x D^2$ of $C_i$ with
$S^1\x({\text{pt on $\bd D^2$}})$ a longitude of $C_i$, {\it i.e.}\/
nullhomologous in
$S^3\setminus C_i$. For any $A\in GL(2;\Z)$ with $\det A=-1$, we the get a
3-manifold
\[ s(L;A)=(S^3\setminus {\text{int}} (N_1\cup N_2))/ A \]
called a {\it{sewn-up link exterior}} by identifying $\bd N_1$ with $\bd N_2$ via a
diffeomorphism inducing $A$ in homology. For
$n\in\Z$, let
$A_n=\bigl(\begin{smallmatrix} -1&0\\-n&1 \end{smallmatrix}\bigr).$
A simple calculation shows that $H_1(s(L;A_n);\Z)=\Z\oplus\Z_{n-2\l}$ where
$\l\,$  is
the linking number in $S^3$ of the two components $C_1$, $C_2$, of $L$.
(See \cite{WB}.)
The second summand is generated by the meridian to either component.

J. Hoste \cite{Hoste} has given a recipe for producing Kirby calculus
diagrams for
$s(L;A_n)$. First we review the notion of a `band sum' on an oriented knot
or link $L$.
Consider a portion of $L$ consisting of a pair of strands, oriented in opposite
directions, and separated by a band $B$. We identify $B$ with an embedding
of $I^2=[0,1]\x
[0,1]$ in $S^3$ such that $B\cap L =(\{0\}\x I)\cup -(\{1\}\x I)$. Let $K$ be
the (oriented) knot or link obtained by trading the segments $B\cap L=
\{0,1\}\x I$ of $\bd
B$ for the complementary oriented segments  $I\x\{0\}$ and $-(I\x\{1\})$.
The process of
exchanging $L$ for $K$ in this fashion is called a {\it{band sum}}.
Associated with the band move are two unknots: $U_{\nu}$, an
unknotted circle which bounds a disk whose interior meets $B=I^2$ in the arc
$\{\frac12\}\x I$ and is disjoint from $L$, and $U_o$, which spans a disk whose
interior meets $B$ in $I\x\{\frac12\}$ and is disjoint from $K$. (See Figure 2.)

\centerline{\unitlength 1cm
\begin{picture}(8,3.75)
\put (1,2){\oval(1.5,1)[l]}
\put (1.2,1){\line(0,1){.4}}
\put (1.2,1.6){\vector(0,1){1.65}}
\put (1.9,1.6){\line(0,1){1.65}}
\put (1.9,1.4){\vector(0,-1){.4}}
\put (2,2){\oval(1.5,1)[r]}
\put (1,1.49){\line(1,0){1}}
\put (1.35,2.49){\line(1,0){.4}}
\multiput (1.2,1.8)(0,.05){10}{\line(1,0){.7}}
\put (2.5,1.25){\small{$U_o$}}
\put (1.4,.6){\small{$L$}}
\put (.6,3){\small{$C_1$}}
\put (2.1,3){\small{$C_2$}}
\put (.8,1.9){\small{$B$}}
\put (3.65,2.05){$\longrightarrow$}
\put (5.2,2.3){\vector(0,1){1}}
\put (5.2,2.3){\line(1,0){.15}}
\put (5.2,1.8){\line(1,0){.15}}
\put (5.2,1.8){\line(0,-1){1}}
\put (5.85,2.1){\oval(.75,1.75)[l]}
\put (5.6,2.3){\line(1,0){.45}}
\put (5.6,1.8){\line(1,0){.45}}
\put (6.05,2.3){\line(0,1){1}}
\put (6.05,1.8){\vector(0,-1){1}}
\put (6.2,2.1){\oval(.75,1.75)[r]}
\put (6.7,2.6){\small{$U_{\nu}$}}
\put (5.45,.6){\small{$K$}}
\put (3.275,0){Figure 2}
\end{picture}}

\begin{prop}[Hoste \cite{Hoste}]\label{JH} Let $L=C_1\cup C_2$ be an
oriented link in $S^3$.
Consider a portion of $L$ consisting of a pair of strands, one from each
component, oriented in opposite directions, and separated by a band $B$.
The band
sum of $C_1$ and
$C_2$ is a knot $K$, and $U_{\nu}$ links $K$ twice geometrically and $0$ times
algebraically. The sewn-up link exterior $s(L;A_n)$ is obtained from
surgery on the the two
component link $K\cup U_{\nu}$ in $S^3$ with surgery coefficient $0$ on
$U_{\nu}$ and
$n-2\l$ on $K$, where $\l$ is the linking number of $C_1$ and $C_2$.
\end{prop}

Next consider a related situation. Let $Z=s(L;A_n)$ where $n=2\l$; so $Z$ has
$H_1(Z;\Z)=\Z\oplus\Z$. Suppose that $B$ is a band in $S^3$ meeting $L$ as in
Proposition~\ref{JH}, with the circle $U_o$, which links $L$ twice
geometrically and
$0$ times algebraically. Then $U_o$ gives rise to a loop, $\bar{U}_o$ in $Z$.
To get a Kirby calculus picture of this situation, apply Hoste's formula.
We obtain a two
component framed link $K\cup U_{\nu}$ with $0$-framing on each. In $S^3$,
$U_o$ bounds a
disk which is disjoint from $K$ and meets a disk spanning $U_{\nu}$ in two
points with
opposite orientations. Thus
$U_o$ bounds a punctured torus in $S^3\setminus (K\cup U_{\nu})$; and so
$\bar{U}_o$
is nullhomologous in $Z$. This means that $\bar{U}_o$ has a naturally defined
longitude in $Z$; so $p/q$- Dehn surgery on $\bar{U}_o$ is well-defined.
Let $Z_0$ denote the result of $0$-surgery on $\bar{U}_o$ in $Z$.

In $Z$, let $m$ denote the meridian circle to the sewn-up torus. Then in
$Z\x S^1$ we
have the torus $T_m=m\x S^1$ of self-intersection $0$. Form the fiber sum
$X\#_T(Z\x S^1)$ of $X$ with $Z\x S^1$ by identifying $T_m$ with the torus
$T$ of
$X$. (In the name of brevity, we shall sometimes make, as in this case, a
mild change
in our notation for fiber sum.) Similarly we can form the fiber sum
$X\#_T(Z_0\x S^1)$, which is seen to be the
result of performing a $(0/1)$-log transform on the torus $\bar{U}_o\x S^1$ in
$X\#_T(Z\x S^1)$. We wish to compute the Seiberg-Witten invariant of this
4-manifold.

By cutting open $s(L;A_n)$ along one of the two components of $L$, we
obtain the exterior
of the link $L$ in $S^3$, and if this is done after performing $0$-surgery on
$U_o$, we obtain a link exterior in $S^2\x S^1$. Now perform the corresponding
task in $X\#_T(Z_0\x S^1)$, removing a torus $C_i\x S^1$. We obtain
$X\#_T(S^2\x S^1\x S^1)$ with a pair of tubular neighborhoods of
self-intersection
$0$ tori removed. Call this manifold $Q(L)$. To reiterate --- $Q(L)$ is
obtained by
performing $0$-surgery on $U_o$ in $S^3\setminus L$, crossing with $S^1$,
and then
fiber-summing with $X$ along $T_m$.  The boundary components of $Q(L)$ have
a natural
framing $\lam$, $\mu$, $\s$, coming from the longitude and meridian of the link
components in $S^3$ and the $S^1$ in the last coordinate. We can re-obtain
$X\#_T(Z_0\x S^1)$ by sewing up the boundary $3$-dimensional tori of
$Q(L)$ using the
matrix $A_n\oplus(1)$.  Instead, let us fill in each of the boundary
components of $Q(L)$
with a copy of $S^1\x D^2\x S^1$. This can be done in many ways. We wish to
do it so that,
using the framings obtained from our copies of $S^1\x D^2\x S^1$, we obtain
$X\#_T(Z_0\x S^1)$ by sewing up the boundary of the resultant manifold with
a neighborhood of
the (new) link $(S^1\x{\text{pt}}\x S^1) \cup (S^1\x{\text{pt}}\x S^1)$
removed using the
matrix
$\pm A_0\oplus(1)$. Using the obvious framing for
$S^1\x D^2\x S^1$, we claim that this is done by gluing each $S^1\x D^2\x
S^1$ to a
component of $\bd Q(L)$ by a diffeomorphism with matrix $B_{\l}\oplus(1)$ where
$B_{\l}= \begin{pmatrix} 0&1\\1&\l \end{pmatrix}$. We shall denote by
$W(L)$ the manifold formed using $B_{\l}\oplus(1)$ to sew in the
neighborhoods of
the tori $C_i' \times S^1 = S^1 \times pt \times S^1.$

Let $V(L)$ be the result of filling in the exterior of $L$ in
$S^3$ via the diffeomorphism $B_{\l}$ on each component. This is just the
result of
$\l$-framed surgery on each component of $L$. Now if we sew up the link
complement
$V(L)\setminus L$ via $-A_0$, we get the result of sewing up $S^3\setminus
L$ using the
diffeomorphism $B_{\l}(-A_0)B_{\l}^{-1} = A_{2\l}.$ Thus $s(L;A_n)=s(L';-A_0)$
where $L'$ is the link $C_1',C_2'.$ Denote by $V_0(L)$ the result of performing
$0$-surgery on $U_o$ in $V(L)\setminus L' = S^3 \setminus L.$ Then
$W(L) = X\#_T (V_0(L) \times S^1)$ and $X\#_T (Z_0 \times S^1) = W(L)_{T_1=T_2}$
where $T_i = C_i' \times S^1.$

\begin{prop}\label{skein2} Suppose that $T$ is a c-embedded torus in $X$. Then
with the above notation,
\[ \sw_{X\#_T(Z_0\x S^1)} = (t^{1/2}-t^{-1/2})^2\cdot\sw_{X_K}\]
as Laurent polynomials, where $K$ is the band sum of $C_1$ and $C_2$ using
the band $B$,
and $t=\exp(2[T])$.
\end{prop}
\begin{proof} Since the matrix $-A_0\oplus(1)$ identifies the tori $C_i'\x S^1$
in the boundary components of $Q(L)$, Theorem~\ref{gluing1} tells us that
$\sw_{X\#_T(Z_0\x S^1)}= \sw_{W(L)_{T_1=T_2}}$ is obtained from the
relative invariant
$\sw_{(W(L);C_1',C_2')}$ by identifying the homology classes in $Q(L)$
represented by the tori $T_i = C_i'\x S^1$.

The gluing diffeomorphism $B_{\l}$ identifies the homology class of a
longitude of $C_i'$
with the  meridian of $C_i$ in $S^3\setminus L$. Thus Theorem~\ref{gluing1}
implies that
\[ \sw_{X\#_T(Z_0\x S^1)}=\sw_{W(L)}\cdot (\t-\t^{-1})^2,\]
where $\t=\exp([m\x S^1])$.

It remains to identify the manifold $W(L)$ as $X_K$. By construction,
$W(L)$ is obtained
from the 3-component link $C_1\cup C_2\cup U_o$ in $S^3$ by performing
$\l$-framed surgery
on $C_1$ and $C_2$ and $0$-framed surgery on $U_o$, crossing with $S^1$,
and fiber summing
to $X$ along $T_m$ and $T$. The result of framed surgery on the 3-component
link is,
by sliding
$C_1$ over $C_2$, seen to be the same as $0$-framed surgery on $K$. Thus
$W(L)=X_K$, and
the handle slide carries the meridian $m$ to a meridian of $K$. Letting
$t=\exp(2[T])=\t^2$, we get
the calculation as claimed.
\end{proof}

\section{The proof of Theorem~\ref{knotthm}\label{proof}}

We first recall a standard technique for calculating the (symmetrized) Alexander
polynomial of a knot.
This uses the skein relation
\begin{equation} \DD_{K_+}(t)=\DD_{K_-}(t)+(t^{1/2}-t^{-1/2})\cdot\DD_{K_0}(t)
\label{boogie}\end{equation}
where $K_+$ is an oriented knot or link, $K_-$ is the result of
changing a single oriented positive (right-handed) crossing in $K_+$ to a
negative
(left-handed) crossing, and $K_0$ is the result of resolving the crossing
as shown in
Figure 3.

Note that if $K_+$ is a knot, the so is $K_-$, and $K_0$ is a 2-component
link.  If $K_+$
is a 2-component link, then so is $K_-$, and $K_0$ is a knot.

\centerline{\unitlength 1cm
\begin{picture}(9,4)
\put (1,1){\vector(1,2){1.12}}
\put (2,1){\line(-1,2){0.45}}
\put (1.45,2.1){\vector(-1,2){.575}}
\put (1.25,0.5){$K_+$}
\put (4,1){\line(1,2){0.45}}
\put (4.55,2.1){\vector(1,2){.575}}
\put (5,1){\vector(-1,2){1.12}}
\put (4.25,0.5){$K_-$}
\put (7,1){\line(1,4){.255}}
\put (8,1){\line(-1,4){.255}}
\put (7.255,2){\vector(-1,4){.3}}
\put (7.745,2){\vector(1,4){.3}}
\put (7.25,0.5){$K_0$}
\put (3.8,-.15){Figure 3}
\end{picture}}

\medskip

The point of using \eqref{boogie} to calculate $\DD_K$ is that $K$ can be
simplified to an
unknot via a sequence of crossing changes of a projection of the oriented
knot or link to
the plane. To describe this well-known technique, consider such a
projection and choose a
basepoint on each component. In the case of a link, order the components.
Say that such a
projection is {\it{descending}}, if starting at the basepoint of the first
component and
traveling along the component, then from the basepoint of the second
component and
traveling along it, {\it{etc.}}, the first time that each crossing is met, it is
met by an
overcrossing. Clearly a link with a descending projection is an unlinked
collection of
unknots. Our goal is to start with a knot $K$ and perform skein moves so as
to build a
tree starting from $K$ and at each stage adding the bifurcation of Figure 4,
where each $K_+$, $K_-$, $K_0$ is a knot or 2-component link, and
so that at the
bottom of the tree, we obtain only unknots, and split links. Then, because
for an unknot
$U$ we have $\DD_U(t)=1$, and for a split link $S$ (of more than one
component) we have
$\DD_S(t)=0$, we can work backwards using \eqref{boogie} to calculate
$\DD_K(t)$.

\centerline{\unitlength .5cm
\begin{picture}(5.5,5.5)
\put (2.25,4){\line(-3,-4){1.5}}
\put (2.75,4){\line(3,-4){1.5}}
\put (2.25,4.3){$K_+$}
\put (0.25,1.2){$K_-$}
\put (4,1.2){$K_0$}
\put (1.3,0){Figure 4}
\end{picture}}

The recipe for constructing the tree is, in the case of a knot, to change
the crossing of
the first `bad' crossing encountered on the traverse described above. In
this case, the
result of changing the crossing is still a knot, and the result of
resolving the crossing
is a 2-component link. In the case of a 2-component link, one changes the
first `bad'
crossing between the two components which is encountered on the traverse.
The result of
changing the crossing is still a 2-component link, and the result of
resolving the
crossing is a knot. In this way we obtain a tree whose top vertex is the
given knot $K$
and which branches downward as in the figure above. We shall call this tree a
{\it{resolution tree}}\/ for the knot $K$. We claim that the tree can be
extended until
each bottom vertex represents an unknot or a split link.

For the projection of an oriented, based knot $K$, let $c(K)$ be the
number of crossings and $b(K)$ be the number of bad crossings encountered
on a traverse
starting at the basepoint. The {\it{complexity}}\/ of the projection is
defined to be
the ordered pair $(c(K), b(K))$. For the projection of an oriented, based
2-component
link $L$, let $c(L)$ be the total number of crossings and let $b(L)$ be the
number of bad
crossings between the two components. Again the complexity is defined to be
$(c(L), b(L))$.
Consider a vertex which represents a knot or 2-component link $A$. Note that
$c(A_-)=c(A_+)$, $b(A_-)<b(A_+)$, and $c(A_0)<c(A_+)$. Thus in the
lexicographic ordering,
$(c(A_-),b(A_-))<(c(A_+),b(A_+))$ and $(c(A_0),b(A_0))<(c(A_+),b(A_+))$.
Now a knot
$K_1$ with $c(K_1)=0$ or with $b(K_1)=0$ is the unknot, and a link $L$ with
$c(L)=0$ is
the unlink and with $b(L)=0$, it is at least a split link. This completes
the proof that
we can construct the resolution tree as described. (We remark that for the
sake of
simplicity we have considered only the case where we have changed a
positive to a negative
crossing in the skein move. Of course, we may as well have to change a
negative to
a positive crossing in order to lower $b$ at various steps, but this does
not change the
proof.)

\medskip

Consider an oriented, based knot $K$ in $S^3$ and a knot projection of $K$.
We shall use the
resolution tree for this projection as a guide for simplifying $X_K$ in a
way which leads
to a calculation of $X_K$. Let us consider the first step, say
$K=K_+\to\{ K_-,\ K_0\}$.
At the crossing of $K_-$ that is in question, there is an unknotted circle
$U$ linking $K$
algebraically $0$ times, so that the result of $+1$ surgery on $U$ turns
$K_-$ into
$K_+=K$. (See Figure 5.)
In $X_{K_-}$ we have the nullhomologous torus $S^1\x U$. (It
bounds the
product of $S^1$ with a punctured torus in $S^3\setminus K_-$.) The fact
that $+1$
surgery on $U$ turns $K_-$ into $K_+$ means that $X_{K_+}$ is the result
of a $(1/1)$-log
transform on $X_{K_-}$ along $S^1\x U$. Let $X_{K_-}(0/1)$ denote the
result of performing
a $(0/1)$-log transform on $S^1\x U$ in $X_{K_-}$. We now use \eqref{log} and
\eqref{log+} to compute $\sw(X_K)$. Two tori are central to this
calculation. Letting
$m_U$  be a meridional circle to $U$, we get the torus $S^1\x m_U$ which
represents the
homology class $\t$ of \eqref{log}. Also, in $X_{K_-}(0/1)$, the boundary
of the punctured torus described above is spanned by a disk, and we obtain
a torus of
self-intersection $0$ representing a class $\s$ such that $\s\cdot\t=1$.
Note that $H_2(X_{K_-}(0/1))=H_2(X_{K_-})\oplus H(\s,\t)$; so \eqref{log+}
applies.
Hence
$ \sw_{X_K}= \sw_{X_{K_-}} + \sw_{X_{K_-}(0/1)}$.

\centerline{\unitlength 1cm
\begin{picture}(3.25,3.75)
\put (1,2){\oval(1.5,1)[l]}
\put (1,1){\line(1,2){.2}}
\put (1.3,1.6){\vector(1,2){.8}}
\put (1.8,1.4){\vector(1,-2){0.2}}
\put (1.45,2.1){\line(-1,2){.575}}
\put (2,2){\oval(1.5,1)[r]}
\put (1,1.49){\line(1,0){1}}
\put (1.35,2.49){\line(1,0){.25}}
\put (1.712,1.575){\line(-1,2){.18}}
\put (2.5,1.25){\small{+1}}
\put (1.35,.6){\small{$K_-$}}
\put (.9,0){Figure 5}
\end{picture}}

\smallskip

Recall that $X_{K_-}(0/1)$ is obtained by performing $0$-framed surgery on
both components
of the link $K_-\cup U$, crossing with $S^1$ and fiber-summing with $X$,
using the torus
$T$ obtained from a meridian of $K_-$ crossed with $S^1$.
Hoste's recipe, \eqref{JH}, allows us to interpret the result of $0$-framed
surgery on
both components of $K_-\cup U$ in $S^3$ as $s(K_0;A_{2\l})$ where $K_0$ is
the 2-component
link obtained by resolving the crossing under consideration (see Figure~6),
and $\l$ is the
linking number of the two components of $K_0$.

\centerline{\unitlength 1cm
\begin{picture}(8,3.75)
\put (1,2){\oval(1.5,1)[l]}
\put (1,1){\line(1,2){.2}}
\put (1.3,1.6){\vector(1,2){.8}}
\put (1.8,1.4){\vector(1,-2){0.2}}
\put (1.45,2.1){\line(-1,2){.575}}
\put (2,2){\oval(1.5,1)[r]}
\put (1,1.49){\line(1,0){1}}
\put (1.35,2.49){\line(1,0){.25}}
\put (1.712,1.575){\line(-1,2){.18}}
\put (2.5,1.25){\small{0}}
\put (1.3,.6){\small{$K_-$}}
\put (.5,3){\small{0}}
\put (6.25,1.7){\vector(1,2){.75}}
\put (6.75,1.7){\line(-1,2){.2}}
\put (6.4,2.25){\line(-1,2){.5}}
\put (6.25,1.5){\line(-1,-2){.35}}
\put (6.25,1.7){\line(1,0){.5}}
\put (6.25,1.5){\line(1,0){.5}}
\put (6.75,1.5){\vector(1,-2){.35}}
\put (3.5,2){$\longrightarrow$}
\put (4.75,2){${\small{s}}\ \Biggl($}
\put (7.8,2){$\Biggr)$}
\put (7.2,2) {\small{$; A_{2\l}$}}
\put (6.3,.6){\small{$K_0$}}
\put (3.3,0){Figure 6}
\end{picture}}

\smallskip

Let $X(s(K_0;A_{2\l}))=(s(K_0;A_{2\l})\x S^1)\#_T X$ where $T$ is the product
of $S^1$ with a
meridian to either component of $K_0$. Because $A_{2\l}$ sends meridians to
meridians,
this definition does not depend on the choice of component. Then
$X_{K_-}(0/1)\cong
X(s(K_0;A_{2\l}))$; so
\begin{equation}
 \sw_{X_{K_+}}=\sw_{X_{K_-}}+\sw_{X(s(K_0;A_{2\l}))},
\label{knot}\end{equation}
mimicking the skein move which gives the second tier of the resolution tree.

Now consider the next stage of the resolution tree and the skein move
$L\to \{L_-,L_0\}$ where $K_0=L=L_+$. This move corresponds to changing a
bad crossing
involving both components of $L$. If the bad crossing under consideration
is, say,
right-handed, then $L=L_+$ can be obtained from $L_-=C_1\cup C_2$ by
$+1$-surgery on an
unknotted circle $U_o$ as in Figure~7.
This means that $X(s(L;A_{2\l}))$ is the result of a $(1/1)$-log
transform on the
torus $S^1\x U_o$ in $X(s(L_-;A_{2\l_-}))$ where $\l_-$ is the linking
number of the
components $C_1$, $C_2$ of $L_-$ and is determined by the fact that
$H_1(s(L_-;A_{2\l_-});\Z)=H_1(s(L;A_{2\l});\Z)=\Z\oplus\Z$.

\centerline{\unitlength 1cm
\begin{picture}(3.25,3.75)
\put (1,2){\oval(1.5,1)[l]}
\put (1,1){\line(1,2){.2}}
\put (1.3,1.6){\line(1,2){.18}}
\put (1.8,1.4){\vector(1,-2){0.2}}
\put (1.55,2.1){\vector(1,2){.55}}
\put (2,2){\oval(1.5,1)[r]}
\put (1,1.49){\line(1,0){1}}
\put (1.35,2.49){\line(1,0){.25}}
\put (1.712,1.575){\line(-1,2){.8}}
\put (2.5,1.25){\small{+1}}
\put (1.3,.6){\small{$L_-$}}
\put (.425,3){\small{$C_1$}}
\put (2.25,3){\small{$C_2$}}
\put (.9,0){Figure 7}
\end{picture}}

\smallskip

In $X(s(L_-;A_{2\l_-}))$, the torus  $S^1\x U_o$ is nullhomologous. This is
precisely the
situation of \eqref{skein2} where $Z=s(L_-;A_{2\l_-})$. We wish to apply
\eqref{log+} to
this situation. Let $X(s(L_-;A_{2\l_-}))(0/1)$ denote the 4-manifold
obtained by performing a
$(0/1)$-log transform on $S^1\x U_o$ in $X(s(L_-;A_{2\l_-}))$. This is the
manifold
$X\#_T(Z_0\x S^1)$ of \eqref{skein2}, and the 3-manifold $Z_0$ is the result of
$0$-surgery on $U_o$ in $s(L_-;A_{2\l_-})$. Let $m_{U_o}$ be a meridional
circle to
$U_o\subset s(L_-;A_{2\l_-})$. The torus \ $m_{U_o}\x S^1$ in
$X(s(L_-;A_{2\l_-}))(0/1)$ is
the $T_0$ mentioned in the statement of \eqref{log}. As we argued in the
proof of
\eqref{skein2}, in $Z$, the loop $U_o$ bounds a punctured torus, and this
gets completed to
a torus of self-intersection $0$ in $Z_0$. Let $\s$ be its homology class in
$X(s(L_-;A_{2\l_-}))(0/1)=X\#_T(Z_0\x S^1)$. Since the class $\s$ satisfies
the hypothesis
of \eqref{log+}, we have
\[ \sw_{X(s(L_+;A_{2\l}))}=\sw_{X(s(L_-;A_{2\l_-}))}+\sw_{X\#_T(Z_0\x S^1)}.\]
Applying \eqref{skein2}, this becomes:
\begin{equation}
\sw_{X(s(L_+;A_{2\l}))}=\sw_{X(s(L_-;A_{2\l_-}))}+(t^{1/2}-t^{-1/2})^2\cdot\sw
_{X_{L_0}}
\label{link}\end{equation}
where $L_0$ is the result of resolving the crossing of $L$ which is under
consideration,
and $t=\exp(2[T])$.

In order to see that this process calculates $\DD_K(t)$, for fixed $X$, we
define a formal Laurent series $\Th$,  which is an invariant of
knots and 2-component links. For a knot $K$, define $\Th_K$ to be the quotient,
$\Th_K=\sw_{X_K}/\sw_X$, and for a
2-component link with linking number $\l$ between its components,
$\Th_L=(t^{1/2}-t^{-1/2})^{-1}\cdot\sw_{X(s(L;A_{2\l}))}/\sw_X$, where as usual
$t=\exp(2[T])$.
It follows from \eqref{knot} and \eqref{link} that for knots or 2-component
links,
$\Th$ satisfies the skein relation
\[ \Th_{K_+} = \Th_{K_-} + (t^{1/2}-t^{-1/2})\cdot \Th_{K_0}. \]
For a split 2-component link, $L$, the 3-manifold $s(L;A_{2\l})$ contains an
essential
2-sphere (coming from the 2-sphere in $S^3$ which splits the link). This
means that
$X(s(L;A_{2\l}))$ contains an essential 2-sphere of self-intersection $0$,
and this
implies that $\sw_{X(s(L;A_{2\l}))}=0$ (see \cite{FS}). Thus for a split
link, $\Th_L=0$.
For the unknot $U$, the manifold $X_U$ is just $X\#_T (S^2\x T^2)=X$, and so
$\Th_U=1$.
Subject to these initial values, the resolution tree and the skein relation
\eqref{boogie} determine
$\DD_K(t)$ for
any knot $K$. It follows that $\Th_K$ is a Laurent polynomial in a single

variable $t$, and $\Th_K(t)=\DD_K(t)$, completing the
proof of Theorem~\ref{K3thm}.

\section{The proof of Theorem~\ref{linkthm}\label{linkproof}}

We first review the axioms which determine the Alexander polynomial of a
link. The  reference for this material is \cite{Turaev}. The fact, proved in
\cite{Turaev}, that we shall use here is that there is but one map $\DN$
which assigns
to each
$n-$component ordered link $L$ in $S^3$ an element of the field
${\mathbf{Q}}(t_1,\dots,t_n)$ with the following properties:

\begin{enumerate}
\item $\DN(L)$ is unchanged under ambient isotopy of the link $L$.
\item If $L$ is the unknot, then $\DN(L)=1/(t-t^{-1})$.
\item If $n \ge 2$, then $\DN(L)\in\Z[t_1,t_1^{-1},\dots,t_n,t_n^{-1}]$.
\item The one-variable function $\tilde{\DN}(L)(t)=\DN(L)(t,t,\dots,t)$ is
unchanged
by a renumbering of the components of $L$.
\item (Conway Axiom). If $L_+$, $L_-$, and $L_0$ are links coinciding (except
possible for the numbering of the components) outside a ball, and inside this
ball have the form depicted in Figure 8, then
\[
\tilde{\DN}(L_+)=\tilde{\DN}(L_-)+(t-t^{-1})\cdot\tilde{\DN}(L_0).
\]
\item (Doubling Axiom).  If the link $L'$ is obtained from the link
$L=\{K_1,\dots,K_n\}$ by replacing the $K_j$ by its $(2,1)-$cable, then
\[\DN(L')(t_1,\dots,t_n)=(T+T^{-1})
\cdot\DN(L)(t_1,\dots,t_{j-1},t_j^2,t_{j+1},\dots,t_n)\] where
$T=t_j\prod_{i\ne j}t_i^{\ell k(K_j,K_i)}$.
\end{enumerate}

\centerline{\unitlength 1cm
\begin{picture}(9,4)
\put (1,1){\vector(1,2){1.12}}
\put (2,1){\line(-1,2){0.45}}
\put (1.45,2.1){\vector(-1,2){.575}}
\put (1.25,0.5){$L_+$}
\put (4,1){\line(1,2){0.45}}
\put (4.55,2.1){\vector(1,2){.575}}
\put (5,1){\vector(-1,2){1.12}}
\put (4.25,0.5){$L_-$}
\put (7,1){\line(1,4){.255}}
\put (8,1){\line(-1,4){.255}}
\put (7.255,2){\vector(-1,4){.3}}
\put (7.745,2){\vector(1,4){.3}}
\put (7.25,0.5){$L_0$}
\put (3.8,-.15){Figure 8}
\end{picture}}
\medskip

The Alexander polynomial $\DD_L$ in Theorem~\ref{linkthm} is just
\[\DD_L(t_1^2,\dots,t^2_n)\doteq\DN(L)(t_1,\dots,t_n) \]
where the symbol $\doteq$ denotes equality up to multiplication by $-1$ and
powers of the variables.

Recall the construction of $E(1)_L$. If $L=\{K_1,\dots,K_n\}$ is an
($n\ge 2$)-component ordered link in $S^3$ and $(\ell_j,m_j)$ denotes the
standard
longitude-meridian pair for the knot $K_j$, we let \[\alpha_L:\pi_1(S^3\setminus
L)\to \Z \]   denote the homomorphism characterized by the property that
$\alpha_L(m_j)=1$ for each
$j=1,\dots,n$. Define  $M_L$ to be the 3-manifold
obtained by performing
$\alpha_L(\ell_j)$ surgery on each component $K_j$ of L. Then,  in
$M_L\x S^1$ we use the smooth tori $T_{m_j}=m_j\x S^1$ to construct
the $n-$fold fiber sum
\[ E(1)_L=(M_L\x S^1)\#_{T_{m_j}=F}\coprod\limits_{j=1}^n E(1),\]
the fiber sum being  performed using the natural
 framings $T_{m_j} \x D^2$ of the neighborhoods of $T_{m_j}=m_j\x S^1$ in $M_j$
and the neighborhoods
$F
\x D^2$ of an elliptic fiber in each copy of $E(1)$ and glued together so
as to preserve
the homology classes
$[{\text{pt}}\x
\bd D^2]$. It is amusing to note that this later condition is unnecessary
in this special
situation. For, since $E(1)\setminus F$ has a big diffeomorphism group, we
can fiber sum
in the $E(1)$ with any gluing map and end up with diffeomorphic $4$-manifolds.

Now let $\DN$ be that function which associates to every ordered
$n$-component link $L$,
the polynomial
$\DN(L)(t_1,\dots,t_n)=\sw_{E(1)_L}(t_1^2,\dots,t_n^2)$, with
$t_j=\exp(2[T_{m_j}])$. We show that $\DN$  satisfies the above stated
axioms. However,
there is a small obstruction to doing this in a straightforward manner:
$\sw_X$ is
only defined when
$b^+>1$ and when
$L$ has but one component $E(1)_L$ has $b^+=1$. Although our
Theorem~\ref{linkthm} is
 stated only for $n\ge 2$, the axioms insist that we consider $1$-component
links. We
overcome this problem as follows. Given an ordered
$n$-component link
$L$ we always fiber sum in the $K3$-surface (rather  than $E(1)$) to $M_L\x
S^1$ along
$T_{m_1}$. In the case that $n\ge 2$, the resulting manifold, which we
temporarily
denote by
$E(2,1)_L$, has
\[\sw_{E(2,1)_L}=(t_1^{1/2}-t_1^{-1/2})\cdot \sw_{E(1)_L}\]
We shall complete the proof of Theorem~\ref{linkthm} by showing that
\[
\DN(L)(t_1,\dots,t_n)=\frac{\sw_{E(2,1)_L}(t_1^2,\dots,t_n^2)}{t_1-t_1^{-1}}
\]
satisfies all the axioms.

Axiom 1 is clear.

For Axiom 2 note that if $L$ is the unknot, then $E(2,1)_L$ is
the $K3$-surface, so that
$\DN(L)(t)=\sw_{E(2,1)_L}(t^2)/(t-t^{-1})=1/(t-t^{-1})$.

To verify Axiom 3, consider an $(n\ge2)$-component link $L$ with components
$K_1,\dots,K_n$.
We need to see that the only possible basic classes of
$E(2,1)_L$ are  the classes $T_{m_i}$ (which are identified in $E(2,1)_L$
with the fiber
classes $F_i$). A Mayer-Vietoris sequence argument shows that
$H_2(E(2,1)_L)\cong{\text{im}}(\varphi)\oplus G$ where
\[ H_2(E(2)\setminus F)\oplus\sum\limits_1^{n-1}H_2(E(1)\setminus F)\oplus
  H_2((S^3\setminus L)\x S^1) \xrightarrow{\varphi} H_2(E(2,1)_L)
  \xrightarrow{\d} G\]
and $G$ is isomorphic to the kernel of \ $\sum\limits_1^nH_1(T^3)\to
H_1((S^3\setminus L)\x
S^1)$. Now
$H_2(E(2)\setminus F)\cong2E_8\oplus 2H\oplus 3(0)$ and
$H_2(E(1)\setminus F)\cong E_8\oplus 3(0)$ where $H$ denotes a hyperbolic
pair. Each copy of
$E_8$ is represented by eight $-2$-spheres in the usual configuration, say
$W$, with
$\bd W=\Sig(2,3,5)$, the Poincar\'e homology sphere. Since $W$ embeds in
$E(2)$, whose only
basic class is $0$, and since $\Sig(2,3,5)$ has positive scalar curvature,
a rudimentary
gluing formula implies that each basic class of $E(2,1)_L$ is orthogonal to
the image of the
$E_8$ summands. Each of the two hyperbolic pairs $H$ is the homology of a
nucleus
in $E(2)\setminus F$ and is generated by  a torus
$\t$ of
self-intersection $0$ and a sphere $\s$ of self-intersection $-2$ which
intersect at one
point. For any basic class $\k$ of $E(2,1)_L$, the adjunction formula
implies $0\ge
\t^2+|\k\cdot\t|$; so $\k$ is orthogonal to $\t$. Also, $\t+\s$ is
represented by another
torus of self-intersection $0$; so $k$ is in fact orthogonal to $H$.
Furthermore, each of
the summands $(0)$ in $H_2(E(1)\setminus F)$  and $H_2(E(2)\setminus F)$ is
represented by a
torus in the boundary of the tubular neighborhood of $F$, and each of these
tori is glued
to a torus in $(S^3\setminus L)\x S^1$ in $E(2,1)_L$. It follows that the
only possible basic classes lying in the image of
$\varphi$ in fact lie in the image of $H_2((S^3\setminus L)\x S^1)$. This
image is spanned
by the classes of the tori $T_{m_i}$, $i=1,\dots,n$ and the tori $V_j$,
$j=1,\dots, n-1$
where $V_i$ is the boundary of a tubular neighborhood of $K_i$ in $S^3$.
The nonzero
elements of $G$ determine classes in $E(2,1)_L$ with nonzero Mayer-Vietoris
boundary. These are
generated by classes $\g_i$, $i=1,\dots,n$ and $\Sig_j$,
$j=1,\dots,n-1$. A representative of $\g_i$ is formed as follows. Let
$S_i$ denote intersection of a Seifert surface for the {\it knot} $K_i$
with the link
exterior. The intersection of $S_i$ with $V_j$ ($j\ne i$) consists of $\l
k(K_i,K_j)$
copies of $m_j$. Each of these is glued to a circle on the fiber
$F$ of the corresponding $E(1)$ or $E(2)$, and this circle bounds in the
elliptic surface. The same is true for the longitude of the knot $K_i$.
The result represents $\g_i$. Note that $\g_i\cdot F_j=\d_{ij}$ and
$\g_i\cdot V_j=0$ for
each $j=1,\dots,n$. The generators $\Sig_j$ are constructed by starting
with an arc $A_j$
in $S^3\setminus L$ from a point on $V_n$ to a point on $V_j$. The boundary
of
$A_j\x S^1$ consists of two circles, and each is identified with a circle in
$H_2(E(1)\setminus F)$ or $H_2(E(2)\setminus F)$. In the elliptic surfaces,
these circles
bound vanishing cycles, disks of self-intersection $-1$. Thus $\Sig_i$ is
represented by a
$-2$-sphere.
We have $\Sig_j\cdot F_i =0$ for all $i$, and $\Sig_j\cdot V_i=\d_{ij}$.
Suppose we have a basic class
\[
\k=\sum\limits_{i=1}^na_iF_i+\sum\limits_{j=1}^{n-1}b_jV_j+\sum\limits_{k=1}
^nc_k\g_k+
    \sum\limits_{\l=1}^{n-1}d_{\l}\Sig_{\l} \]
Since $F_i$ is a torus of self-intersection $0$, the adjunction inequality
implies that
$\k\cdot F_i=0$, {\it{i.e.}} that $c_i=0$. Similarly, $V_j$ is a torus of
self-intersection
$0$; so
$0=\k\cdot V_j=d_j$. Also $\Sig_{\l}+V_{\l}$ is represented by a torus of
self-intersection
$0$; thus $0=\k\cdot (\Sig_{\l}+V_{\l})=\k\cdot\Sig_{\l}=b_{\l}$.

Axiom 4 is clear.

Axiom 5 is verified in the spirit of Theorem~\ref{knotthm}. However, we
must first
construct an auxiliary manifold $\bar{E}(1)_L$ as follows. In $E(2,1)_L$, $L$ an
$n$-component link, let
$F_1,\dots,F_n$ be tori with the torus $F_1$ the elliptic fiber in the
$K3$-surface and, for
$j\ge 1$, $F_j$ the elliptic fiber in the
$(j-1)^{st}$ copy of $E(1)$. (Note that $F_i=T_{m_i}$ in $E(2,1)_L$.) Now
perform $(n-1)$
internal fiber sums, identifying $F_1$ with
$F_2$, a parallel copy of $F_2$ with $F_3$, and so on. The homology classes
represented by the
$F_j$ in $\bar{E}(1)_L$ are all equal, and we denote this homology class by
$[F]$. Let
$t=\exp(2[F])$. It follows from the gluing formula Theorem~\ref{gluing1} that
\[\sw_{\bar{E}(1)_L}(t)=\sw_{E(2,1)_L}(t,\dots,t)\cdot
(t^{1/2}-t^{-1/2})^{(2n-2)},\] or by defining \
$\widetilde{\sw}_{E(1)_L}(t)=\sw_{E(2,1)_L}(t,\dots,t)$,
\[\widetilde{\sw}_{E(1)_L}(t)=
\frac{\sw_{\bar{E}(1)_L}(t)}{(t^{1/2}-t^{-1/2})^{(2n-2)}}.\]

Suppose now that $L_+$, $L_-$, and $L_0$ are links which coincide (except
possibly for the numbering of the components) outside a ball, and inside this
ball have the form depicted in Figure~8. Furthermore, assume $L_\pm$ has $n$
components. Then, as in the proof of Theorem~\ref{knotthm}, there is a
nullhomologous torus
$T$ in
$\bar{E}(1)_{L_-}$ so that
$\bar{E}(1)_{L_+}$ is the result of a
$(1/1)$-log transform on T. By the log transform formula (Theorems~\ref{log} and
\ref{log+}),
\[\sw_{\bar{E}(1)_{L_+}}=\sw_{\bar{E}(1)_{L_-}}+\sw_{\bar{E}(1)_{L_-}(0/1)}.\]

There are two cases. For the first case, the two strands of $L_+$ are from
distinct
components $K_{i_1}$, $K_{i_2}$ of $L_+$. The manifold
$\bar{E}(1)_{L_-}(0/1)$ is obtained
as follows: Perform surgeries on the link components $\{ K_i\}$ of $L_-$
with surgery
coefficient $\a_{L_-}(\l_i)$ on $K_i$ and perform $0$-surgery on the
unknotted component
$U$ which links $K_{i_1}$ and  $K_{i_2}$ as shown in Figure~9.

\centerline{\unitlength 1cm
\begin{picture}(8,4.5)
\put (1,2){\oval(1.5,1)[l]}
\put (1.2,1.6){\line(0,1){1.05}}
\put (1.2,2.85){\oval(.35,.25)[b]}
\put (1.2,2.8){\line(0,1){1}}
\put (1.075,2.9){\oval(.1,.05)[tl]}
\put (1.325,2.9){\oval(.1,.05)[tr]}
\put (1.2,1.35){\line(0,-1){.5}}
\put (2,2){\oval(1.5,1)[r]}
\put (1,1.49){\line(1,0){1}}
\put (1.35,2.49){\line(1,0){.25}}
\put (1.8,1.6){\line(0,1){1.05}}
\put (1.8,1.35){\line(0,-1){.5}}
\put (1.8,2.8){\line(0,1){1}}
\put (1.675,2.9){\oval(.1,.05)[tl]}
\put (1.925,2.9){\oval(.1,.05)[tr]}
\put (1.8,2.85){\oval(.35,.25)[b]}
\put (2.5,1.15){0}
\put (.325,3.5){{\Small{$\a(\l_{i_1})$}}}
\put (1.95,3.5){{\Small{$\a(\l_{i_2})$}}}
\put (.45,2.75){{\Small{$m_{i_1}$}}}
\put (2.1,2.75){{\Small{$m_{i_2}$}}}
\put (.65,.75){{\small{$K_{i_1}$}}}
\put (1.9,.75){{\small{$K_{i_2}$}}}
\put (3.65,2){$\longrightarrow$}
\put (6.3,1.6){\oval(.6,.5)[t]}
\put (6.3,2.4){\oval(.6,.5)[b]}
\put (5.99,1.6){\line(0,-1){.7}}
\put (6.6,1.6){\line(0,-1){.7}}
\put (5.875,2.9){\oval(.1,.05)[tl]}
\put (6.125,2.9){\oval(.1,.05)[tr]}
\put (6,2.85){\oval(.35,.25)[b]}
\put (6.475,2.9){\oval(.1,.05)[tl]}
\put (6.725,2.9){\oval(.1,.05)[tr]}
\put (6.59,2.85){\oval(.35,.25)[b]}
\put (5.99,2.4){\line(0,1){.2}}
\put (6.59,2.4){\line(0,1){.2}}
\put (5.99,2.8){\line(0,1){1}}
\put (6.59,2.8){\line(0,1){1}}
\put (5.3,2.75){{\Small{$m_0$}}}
\put (6.9,2.75){{\Small{$m_0$}}}
\put (4.8,3.5){{\Small{$\a_{L_0}(\l_0)$}}}
\put (6.65,.75){{\small{$K_0$}}}
\put (3.3,.25){Figure 9}
\end{picture}}

 Let $M^3$ be the resulting
3-manifold. Next form
\begin{equation} E(2,1)_{L_-}(0/1)=(M^3\x S^1)
\#_{T_{m_1}=F_1}E(2)\#_{T_{m_2}=F_2}E(1)\#\cdots\#_{T_{m_n}=F_n}E(1).
\label{Ebar0}
\end{equation}
Finally, $\bar{E}(1)_{L_-}(0/1)$ is obtained from $ E(2,1)_{L_-}(0/1)$ by
performing $n-1$
internal fiber sums, as described above. If we slide the handle
corresponding to $K_{i_2}$
over the handle corresponding to $K_{i_1}$ then we obtain a new Kirby
calculus description
of $M^3$: it is obtained from surgery on the link $L_0$ with surgery
coefficients again
given by the $\a_{L_0}(\l_i)$, as in Figure~9. (Note that if the new
component is called
$K_0$ then  for $j\ne i_1$, $i_2$, the linking number
$\l k_{L_0}(K_0,K_j)=\l k_{L_-}(K_{i_1},K_j) + \l k_{L_-}(K_{i_2},K_j)$; so
the total
linking number for the longitude $\l_0$ is
$\a_{L_0}(\l_0)=\a_{L_-}(\l_{i_1})+\a_{L_-}(\l_{i_2})-2\l
k_{L_-}(K_{i_1},K_{i_2})$ which
is exactly the framing which is given to $K_0$ by the handle slide.)

Now the link $L_0$ has $n-1$ components.
We see that the difference between the constructions for
$\bar{E}(1)_{L_-}(0/1)$ and
$\bar{E}(1)_{L_0}$ is that an extra copy of $E(1)$ needs to be fiber-summed into
$\bar{E}(1)_{L_0}$ and then an extra internal fiber sum needs to be
performed on the
result, in order to get $\bar{E}(1)_{L_-}(0/1)$. Theorem~\ref{eiglue}
and Theorem~\ref{gluing1} then imply that
$\sw_{\bar{E}(1)_{L_-}(0/1)}(t)=(t^{1/2}-t^{-1/2})^3\sw_{\bar{E}(1)_{L_0}}(t)$,
where one factor $(t^{1/2}-t^{-1/2})$ comes from the extra fiber sum with
$E(1)$,
and $(t^{1/2}-t^{-1/2})^2$ comes from the extra internal fiber sum. We have:
\begin{eqnarray*}\widetilde{\sw}_{\bar{E}(1)_{L_+}}(t)&=&
\frac{\sw_{\bar{E}(1)_{L_+}}(t)}
{(t^{1/2}-t^{-1/2})^{(2n-2)}}
\\&=&
\frac{ \sw_{\bar{E}(1)_{L_-}}(t)+\sw_{\bar{E}(1)_{L_-}(0/1)}(t) }
{(t^{1/2}-t^{-1/2})^{(2n-2)}}
\\&=& \widetilde{\sw}_{\bar{E}(1)_{L_-}}(t)+
\frac{(t^{1/2}-t^{-1/2})^3\cdot\sw_{\bar{E}(1)_{L_0}}(t)}
{(t^{1/2}-t^{-1/2})^{(2n-2)}}\\&= & \widetilde{\sw}_{\bar{E}(1)_{L_-}}(t)+
\frac{(t^{1/2}-t^{-1/2})\cdot\sw_{\bar{E}(1)_{L_0}}(t)}
{(t^{1/2}-t^{-1/2})^{(2n-4)}}\\& =&
\widetilde{\sw}_{\bar{E}(1)_{L_-}}(t)+(t^{1/2}-t^{-1/2})\cdot
\widetilde{\sw}_{\bar{E}(1)_{L_0}}(t),
\end{eqnarray*} as desired.

For the second case, the two strands of $L_+$ are from the same
component (say the $j^{th}$) of $L_+$; so $L_0$ has $(n+1)$ components.
Then $\bar{E}(1)_{L_-}(0/1)$ is obtained by first performing surgeries on
the components
$\{ K_i\}$ of $L_-$ with surgery coefficient $\a_{L_-}(\l_i)$ on $K_i$ and then
performing $0$-surgery on an unknotted circle which links $K_{j}$ twice
geometrically and $0$-times algebraically as in Figure~10.

\centerline{\unitlength 1cm
\begin{picture}(8,4)
\put (1,2){\oval(1.5,1)[l]}
\put (1,1){\line(1,2){.2}}
\put (1.3,1.6){\line(1,2){.6}}
\put (1.875,3.05){\oval(.1,.05)[tl]}
\put (2.125,3.05){\oval(.1,.05)[tr]}
\put (2,3){\oval(.35,.25)[b]}
\put (1.95,2.95){\vector(1,2){.25}}
\put (2.35,2.9){\small{$m_j$}}
\put (1.8,1.4){\vector(1,-2){0.2}}
\put (1.45,2.1){\line(-1,2){.575}}
\put (2,2){\oval(1.5,1)[r]}
\put (1,1.49){\line(1,0){1}}
\put (1.35,2.49){\line(1,0){.25}}
\put (1.712,1.575){\line(-1,2){.18}}
\put (2.5,1.25){0}
\put (1.3,.6){$K_j$}
\put (0,3){\small{$\a(\l_j)$}}
\put (3.5,2){$\longrightarrow$}
\put (6.25,1.7){\line(1,2){.55}}
\put (6.8,3.05){\oval(.1,.05)[tl]}
\put (7.025,3.05){\oval(.1,.05)[tr]}
\put (6.9,3){\oval(.35,.25)[b]}
\put (6.85,2.95){\vector(1,2){.25}}
\put (6.75,1.7){\line(-1,2){.2}}
\put (6.4,2.25){\line(-1,2){.5}}
\put (6.25,1.5){\line(-1,-2){.35}}
\put (6.25,1.7){\line(1,0){.5}}
\put (6.25,1.5){\line(1,0){.5}}
\put (6.75,1.5){\vector(1,-2){.35}}
\put (4.75,2){$s\ \Biggl($}
\put (7.75,2){$\Biggr)$}
\put (7.2,2) {$; A_n$}
\put (7.25,.6){$K''$}
\put (5.5,3.25){$K'$}
\put (7.25,2.9){\small{$m'$}}
\put (3.2,.25){Figure 10}
\end{picture}}

\noindent This gives a 3-manifold, $M^3$.
Then we form $E(2,1)_{L_-}(0/1)$ as given by \eqref{Ebar0}, and
$\bar{E}(1)_{L_-}(0/1)$ is obtained from this by performing $(n-1)$
internal fiber sums.

Denote the components of $L_0$ by
$K_1,\dots,K_{j-1},K',K'',K_{j+1},\dots,K_n$. Let $L_j$
denote the 2-component link $L_j=\{ K',K''\}$. It follows from Hoste's
theorem \eqref{JH}
that $M^3$ may be obtained from the sewn-up link exterior $s(L_j;A_n)$ by
further surgering
the $K_i$, $i\ne j$ with framing
$\a_{L_-}(\l_i)$  (where $n=\a_{L_-}(\l_j)+2\l k(K',K'')$). Again see Figure~10.
Because $\a_{L_0}(\l')+\a_{L_0}(\l'')=\sum\limits_{i\ne j}\l k(K_i,K_j)+2\l
k(K',K'')=n$,
and also $\a_{L_-}(\l_i)=\a_{L_0}(\l_i)$, $i\ne j$, the discussion preceding
Proposition~\ref{skein2} relates the Seiberg-Witten invariant of
$E(2,1)_{L_-}(0/1)$ with the Seiberg-Witten invariant of $E(2,1)_{L_0}$. We
need to keep
in mind, however, that because $L_0$ has $n+1$ components, there is an
extra fiber sum
with $E(1)$ in the construction for $E(2,1)_{L_0}$. Thus,
\[
\sw_{E(2,1)_{L_-}(0/1)}=(t_j^{1/2}-t_j^{-1/2})\cdot\sw_{E(2,1)_{L_0}}|_{{}_{
t'=t''=t_j}}\]
Furthermore, $\bar{E}(1)_{L_0}$ has one more internal fiber sum than does
$\bar{E}(1)_{L_-}(0/1)$; so
\[
\sw_{\bar{E}(1)_{L_-}(0/1)}(t)=\frac{1}{(t^{1/2}-t^{-1/2})}\cdot\sw_{\bar{E}
(1)_{L_0}}(t).
\]
Thus
\begin{eqnarray*}
\widetilde{\sw}_{\bar{E}(1)_{L_+}}(t)&=&\frac{\sw_{\bar{E}(1)_{L_+}}(t)}
{(t^{1/2}-t^{-1/2})^{(2n-2)}}
\\&=&
\frac{\sw_{\bar{E}(1)_{L_-}}(t)+\sw_{\bar{E}(1)_{L_-}(0/1)}(t)}
{(t^{1/2}-t^{-1/2})^{(2n-2)}}\\
&=&\widetilde{\sw}_{\bar{E}(1)_{L_-}}(t)+
\frac{\frac{1}{(t^{1/2}-t^{-1/2})}\cdot\sw_{\bar{E}(1)_{L_0}}(t)}
{(t^{1/2}-t^{-1/2})^{(2n-2)}}\\&= & \widetilde{\sw}_{\bar{E}(1)_{L_-}}(t)+
\frac{(t^{1/2}-t^{-1/2})\cdot\sw_{\bar{E}(1)_{L_0}}(t)}
{(t^{1/2}-t^{-1/2})^{2n}}\\& =&
\widetilde{\sw}_{\bar{E}(1)_{L_-}}(t)+(t^{1/2}-t^{-1/2})\cdot
\widetilde{\sw}_{\bar{E}(1)_{L_0}}(t),
\end{eqnarray*} completing the proof of Axiom 5.

Finally, to verify Axiom 6 we first note that $E(2,1)_{L'}$ is obtained
from $E(2,1)_{L}$ by performing an order $2$  log transform on the torus
$\bar{T_j}=
\bar{K_j}\x S^1$ where  $\bar{K_j}$ is the core of the surgered $K_j$ in $M_L$.
The homology class of the resulting meridian $m_j^\prime$ is twice that of
$m_j$ so that
$[T_{m_j^\prime}]=2[T_{m_j}]$. The log transform formulas of
\cite{rat}, then state that
\begin{multline*}
\sw_{E(2,1)_{L'}}(t_1,\dots,t_{j-1},t^\prime_j,t_{j+1},\dots,t_n)\\=(\bar{t}
_j^{\ 1/2}+
\bar{t}_j^{\
-1/2})\cdot\sw_{E(2,1)_L}(t_1,\dots,t_{j-1},t^2_j,t_{j+1},\dots,t_n).
\end{multline*}
The result now follows since
\[ [\bar{T_j}]=[T_j]+\sum_{i\ne j}\ell k(K_j,K_i)[T_i].\]

\section{Examples with $b^+=1$ \label{b^+=1}}

In this section we shall discuss examples which have $b^+=1$. For such
manifolds,
the Seiberg-Witten invariant depends on a choice of metric and self-dual
2-form as
follows. Let $X$ be a simply connected oriented 4-manifold with $b_X^+=1$
with a given
orientation of $H^2_+(X;\R)$ and a given metric $g$. Since $b^+_X=1$, there
is a unique
$g$-self-dual harmonic 2-form $\o_g\in H^2_+(X;\R)$ with $\o_g^2=1$ and
corresponding
to the positive orientation. Fix a characteristic cohomology
class $k\in H^2(X;\Z)$. Given a pair $(A,\psi)$, where
$A$ is a connection in the complex line bundle corresponding to
$k$ and $\psi$ a section of the bundle $W^+$ of self-dual spinors for
the associated $spin^{\,c}$ structure, the perturbed Seiberg-Witten
equations are:
\begin{gather}
D_A\psi = 0 \\
F_A^+  = q(\psi)+i\eta^+ \notag\label{SWeqn}
\end{gather}
where $F_A^+$ is the self-dual part of the curvature of $A$ ,
$D_A$ is the twisted Dirac operator, $\eta^+$, is a
self-dual 2-form on $X$, and
$q$ is a quadratic function. Write $\sw_{X,g,\eta^+}(k)$ for the
corresponding invariant evaluated on the class $k$
($=\frac{i}{2\pi}[F_A]$). As the pair
$(g,\eta^+)$ varies, $\sw_{X,g,\eta^+}(k)$ can change only at those pairs
$(g,\eta^+)$ for which there are solutions of \eqref{SWeqn} with $\psi=0$.
These
solutions occur for pairs $(g,\eta^+)$ satisfying $(2\pi k+\eta^+)\cdot\o_g=0$.
This last equation defines a codimension 1 subspace (`wall') in
$H^2(X;\R)$. The point
$\o_g$ lives in the double cone $\C_X=\{\a\in H^2(X;\R)|\a\cdot\a>0\}$, and, if
$(2\pi k+\eta^+)\cdot\o_g\ne 0$ for a generic $\eta^+$,
$\,\sw_{X,g,\eta^+}(k)$ is
well-defined, and its value depends only on the sign of $(2\pi
k+\eta^+)\cdot\o_g$.
A useful lemma, which follows from the Cauchy-Schwarz inequality (see
\cite{LL}) is:
\begin{lem}\label{lightcone} Suppose that $\a$ and $\b$ are nonzero elements of
$H^2(X;\R)$ which lie in the closure of the same component of $C_X$. Then
$\a\cdot\b\ge0$
with equality if and only if $\a=\lam\b$ for some $\lam>0$.
\end{lem}
It follows from this lemma that that $\sw_{X,g,\eta^+}(k)$  depends only on the
component of $\C_X$ which contains the $g$-self-dual projection of $2\pi
k+\eta^+$.
Furthermore, if the
$g_i$-self-dual projections of $2\pi k+\eta_i^+$ lie in different components of
$\C_X$ and are oriented so that
$(2\pi k+\eta_1^+)\cdot\omega_{g_1}>0>(2\pi k + \eta_0^+)\cdot\omega_{g_0}$,
then
\begin{equation}\label{WC}
\sw_{X,g_1,\eta_1^+}(k)-\sw_{X,g_0,\eta_0^+}(k)=(-1)^{\frac12\d(k)+1}
\end{equation}
where $\d(k)=\frac14(k^2-(3\text{sign}+2\text{e})(X))$ is the formal
dimension of the
moduli spaces.  Thus, as the pair $(g,\eta^+)$ is varied, there are exactly
two values of
$\sw_{X,g,\eta^+}(k)$. Furthermore, in case $b^-\le 9$, for any fixed $a\in
H_2(X;\Z)$
with $\d(a)=a^2+b^--9\ge0$, the self-dual projections of $2\pi a$
all lie in the same component of $\C_X$; so, if $\,a=k\,$ is
characteristic, then for
small enough perturbations, the Seiberg-Witten invariants agree,
independent of metric
\cite{KM,S2}.

Suppose that $X$ contains a smooth essential torus $T$ of
self-intersection $0$. By Lemma~\ref{lightcone}, the class $[T]$ orients
$\C_X$ by declaring $\C^+_X$ to be the component of $\C_X$ which contains
classes
$\a$ with $\a\cdot [T]>0$. Denote the other component by $\C^-_X$ and the
corresponding Seiberg-Witten invariants by $\sw_X^{\pm}$; {\it{i.e.}}
$\sw_X^+(k) = \sw_{X,g,\eta^+}(k)$ where
the $g$-self-dual projection of $2\pi k+\eta^+$ has positive intersection
with $[T]$, and
we define $\sw^-_X(k)$ similarly. The {\it {$[T]^{\perp}$-restricted
Seiberg-Witten
invariants}} are defined to be
\[\sw_{X,T}^{\pm}=\sum_{k\cdot [T]=0}\sw_X^{\pm}(k)\,\exp(k).\]
When $k^2\ge0$ and $k\cdot[T]=0$, Lemma~\ref{lightcone} implies that
$k=\lam [T]$. In particular, if $b^-_X\le9$ and $\d(k)\ge0$, then $k=\lam
[T]$ if
$k\cdot [T]=0$.

\begin{lem}\label{equal!} Let $X$ be a simply connected smooth 4-manifold with
$b^+_X=1$, and suppose that $X$ contains a smooth homologically nontrivial
torus $T$ of
self-intersection $0$. Let $t=\exp(2[T])$. Then
\[ (t^{1/2}-t^{-1/2})\cdot\sw_{X,T}^+=(t^{1/2}-t^{-1/2})\cdot\sw_{X,T}^-.\]
\end{lem}
\begin{proof} The coefficient of $\exp(k)$ in
$(t^{1/2}-t^{-1/2})\cdot\sw_{X,T}^+$ is
\[ c^+(k)=\sw^+_X(k-[T])-\sw^+_X(k+[T]). \]
Since $k\cdot [T]=0$, we have $(k-[T])^2=(k+[T])^2$; so the
wall-crossing formula \eqref{WC} implies that
$c^+(k)=\sw^-_X(k-[T])-\sw^-_X(k+[T])$, the coefficient of
$\exp(k)$ in $(t^{1/2}-t^{-1/2})\cdot\sw_{X,T}^-$.
\end{proof}

For example, consider the case of the rational elliptic surface
$E(1)$ and its fiber class $[F]$. For a Kahler metric $g$ on $E(1)$, the Kahler
form $\,\o\,$ is self-dual, and since $F$ is a complex curve, $\o\cdot
[F]>0$. So the
self-dual projection of $[F]$ is a positive multiple of $\,\o$, and we see
that the
small-perturbation component of $\C_{E(1)}$ is $\C^+_{E(1)}$ for $n[F]$,
$n>0$, and
similarly is $\C^-_{E(1)}$ for $n[F]$, $n<0$.  Since $E(1)$ carries a
metric of positive
scalar curvature, this means $\sw^+_{E(1)}(n[F])=0$ for $n>0$, and
$\sw^-_{E(1)}(n[F])=0$ for $n<0$.
The wall-crossing formula \eqref{WC} implies that (up to an overall sign)
$\sw^-_{E(1)}((2n+1)[F])-\sw^+_{E(1)}((2n+1)[F])=1$. Hence (up to that same
sign)
\[ \sw^-_{E(1)}((2n+1)[F])=\begin{cases} 1,\ \ &n\ge0\\ 0,& n<0.
\end{cases} \]
In other words,
\[ \sw^-_{E(1),F}=\sum_{n=0}^{\infty}t^{(2n+1)/2} \]
where $t=\exp(2[F])$. Similarly,
\[ \sw^+_{E(1),F}=-\sum_{n=0}^{\infty}t^{-(2n+1)/2}. \]
Notice that in this case,
$(t^{1/2}-t^{-1/2})\cdot\sw_{E(1),F}^+$ and
$(t^{1/2}-t^{-1/2})\cdot\sw_{E(1),F}^-$ both
equal $-1$, as Lemma~\ref{equal!} demands.

Let $X$ be a simply connected smooth 4-manifold with $b^+_X=1$, and suppose
that $X$
contains a c-embedded oriented torus $T$. Suppose further that $\pi_1(X\setminus
T)=1$. Then for any knot $K$ in $S^3$ we can form $X_K$, which is
homeomorphic to $X$.
Since $[T]=[T_m]$ orients $\C_{X_K}$, we have invariants $\sw^{\pm}_X$ and
$\sw^{\pm}_{X_K}$. Our discussion of \S~\ref{proof} applies in this case as
well, once we
have the analogues of the log transform and gluing formulas of
\S~\ref{backgd}. These
formulas are also due to Taubes and Morgan, Mrowka, and Szabo, {\it{cf}}.
\cite{S2}.

\begin{ther} [Morgan, Mrowka, and Szabo \cite{MMS}, Taubes] Consider the
simply connected
4-manifolds
$X_1$ and $X_2$ with
$b^+_{X_1}=1$ and $b^+_{X_2}\ge 1$. Suppose that $X_1$ and $X_2$ contain
c-embedded tori
$T_1$ and $T_2$.  Then
\[ \sw_{X_1\#_{T_1=T_2}X_2}=
\begin{cases}
\sw^{\pm}_{X_1,T_1}\cdot\sw^{\pm}_{X_2,T_2}\cdot(t^{1/2}-t^{-1/2})^2,\ \
&\text{if\ \ }
b^+_{X_2}=1 \\
\sw^{\pm}_{X_1,T_1}\cdot\sw_{X_2}\cdot(t^{1/2}-t^{-1/2})^2,\ \ &\text{if\ \ }
b^+_{X_2}>1 \end{cases}\]
where $t=\exp(2[T])$ and $[T]=[T_1]=[T_2]\in H_2(X_1\#_{T_1=T_2}X_2)$.
\end{ther}

It follows from Lemma~\ref{equal!} that these formulas actually make sense. The
restriction to $\sw^{\pm}_{X_1,T_1}$ is accounted for by the fact that
$b^+_{X_1\#_{T_1=T_2}X_2}\ge 3$, and so any basic class in
$H_2(X_1\#_{T_1=T_2}X_2)$ must
be orthogonal to $[T]$.

In order to state the internal fiber sum formula, let $X$ be a simply connected
4-manifold with $b^+_X=1$. Suppose that $X$ contains a pair of c-embedded tori
$T_1$, $T_2$, and let
$X_{T_1,T_2}$ denote the internal fiber sum. Denote by
$\sw^{\pm}_{X,T_1,T_2}$ the
$([T_1],[T_2])^{\perp}$-restricted Seiberg-Witten invariants of $X$,
\[\sw^{\pm}_{X,T_1,T_2}=\sum_{\substack{k\cdot [T_i]=0\\
i=1,2}}\sw_X^{\pm}(k)\,\exp(k).\]

\begin{ther}[Morgan, Mrowka, and Szabo \cite{MMS}, Taubes] Let $X$ be a
simply connected
4-man\-ifold with $b^+_X=1$. Suppose that $X$ contains a pair of c-embedded tori
$T_1$, $T_2$  representing homology classes $[T_1]$ and $[T_2]$.
\[ \sw_{X_{T_1,T_2}}={(\sw^{\pm}_{X,T_1,T_2})|}_{\t_1=\t_2}\cdot
(t^{1/2}-t^{-1/2})^2 \]
where $\t_i=\exp([T_i])$ and $t=\exp(2[T])$.
\end{ther}

Of course we also need the log transform formula:

\begin{ther}[Morgan, Mrowka, and Szabo \cite{MMS}, Taubes] Let $Y$
be a smooth 4-manifold with $b^+=1$, and suppose that $Y$ contains a
nullhomologous torus
$T$. Let $\t$ be the homology class
of $T_0$ in $Y(0/1)$. For each characteristic homology class $\a\in
H_2(Y;\Z)$ and
$p\ne0$,
\[\sw^{\pm}_{Y(p/q)}(\a)= p\sw^{\pm}_Y(\a) +
q\sum\limits_{i=-\infty}^{\infty}\sw_{Y(0/1)}(\a+2i\t).\]
\end{ther}

Now the arguments of \S~\ref{proof} go through verbatim to give:

\begin{thm}\label{b+1K} Let $X$ be a simply connected smooth 4-manifold with
$b^+_X=1$, and suppose that $X$ contains a c-embedded oriented torus $T$.
Suppose further that $\pi_1(X\setminus T)=1$.
For $t=\exp(2[T])$ the $[T]^{\perp}$-restricted Seiberg-Witten invariants
of  $X_K$ are
\[ \sw^{\pm}_{X_K,T}=\sw^{\pm}_{X,T}\cdot\DD_K(t).\]
\end{thm}

As in the $b^+\ge3$ case, if $X$ is a symplectic $4$-manifold with $b^+=1$
containing a
symplectic embedded torus $T$ of self-intersection $0$ and
$K$ is a fibered knot, then $X_K$ is a symplectic $4$-manifold. Conversely,
if $X$ is symplectic, choose a metric $g$ for $X$ so that the symplectic
form $\o$ is
self-dual, and  let $\k_X$ denote the canonical class. In \cite{T1,T2},
Taubes shows that
for $r<<0$, we have $\sw_{X,g,r\o}(-\k_X)=\pm1$
and also that if  $\sw_{X,g,r\o}(k)\ne0$ then
$-\k_X\cdot\o\le k\cdot\o$, with equality only when $k=-\k_X$.
Now let $b^+_X=1$, and let $T$ be an embedded torus of self-intersection $0$,
such that $[T]\cdot\o>0$. For example, this holds if $T$ is symplectically
embedded. Then for any $k\in H_2(X;\R)$, we have
$(2\pi k^++r\o)\cdot [T]<0$ for $r<<0$. This means that
$\sw_{X,g,r\o}(k)=\sw^-_X(k)$.
Hence:

\begin{prop}[Taubes] \label{CT} Let $X$ be a symplectic $4$-manifold with
$b^+=1$
containing an embedded torus $T$ of self-intersection $0$ such that
$[T]\cdot\o>0$. Then
\[ \sw_X^-(-\k_X)=\pm1 \]
and if $\sw^-_{X}(k)\ne0$, then
\[ -\k_X\cdot\o\le k\cdot\o \]
 with equality only when $k=-\k_X$. \qed
\end{prop}

\begin{lem}\label{almfin}  Let $X$ be a simply connected smooth 4-manifold with
$b^+_X=1$, and suppose that $X$ contains a c-embedded oriented torus $T$
 and that $\pi_1(X\setminus T)=1$.
Suppose also that $\sw_{X\#_{T=F}E(1)}\ne0$. Then there is a Laurent polynomial
$S_X$ such that
\[ \sw^{\pm}_{X,T} = S_X\cdot\sum_{n=0}^{\infty}t^{(2n+1)/2} \ \ {\text{or}}\ \
 S_X\cdot\sum_{n=0}^{\infty}t^{-(2n+1)/2}.\]
\end{lem}
\begin{proof} The gluing formula implies that
\[\sw_{X\#_{T=F}E(1)}=\sw^{\pm}_{X,T}\cdot\sw^{\pm}_{E(1),F}\cdot(t^{1/2}-t^
{-1/2})^2
=\sw^{\pm}_{X,T}\cdot(t^{1/2}-t^{-1/2})\]
by Lemma~\ref{equal!} and the above calculation of $\sw^{\pm}_{E(1),F}$.
Because $X\#_{T=F}E(1)$ has $b^+=3$, its Seiberg-Witten invariant is a Laurent
polynomial, $S_X\ne 0$; so the lemma follows.
\end{proof}

\begin{lem}\label{wcr} Suppose that $b^+_X=1$ and $\a\in H^2(X)$, then for any
pair $(g,\eta^+)$ with $\eta^+\cdot\o_g\ne0$ and $|r|>>0$, we have
\[ |\,\sw_{X,g,r\eta^+}(\a)\pm\sw_{X,g,r\eta^+}(-\a)\,|=1.\]
\end{lem}
\begin{proof} Since $F_{\bar{A}}=-F_A$ and $q(\bar{\psi})=-q(\psi)$, if
$(A,\psi)$ is a
solution to the Seiberg-Witten equations for $\a$ corresponding to
$(g,r\eta^+)$, then
$(\bar{A},\bar{\psi})$ is a solution to the equations for $-\a$
corresponding to $(g,-r\eta^+)$. Thus
\[\sw_{X,g,-r\eta^+}(-\a)=(-1)^{(\text{sign}+\text{e})(X)/4}\sw_{X,g,r\eta^+
}(\a).\]
For $|r|>>0$, the signs of
$(-2\pi\a-r\eta^+)\cdot\o_g$ and $(-2\pi\a+r\eta^+)\cdot\o_g$ are opposite.
Thus the
wall-crossing formula implies that
\[ |\,\sw_{X,g,-r\eta^+}(-\a)\pm\sw_{X,g,r\eta^+}(-\a)\,|=1\]
and the lemma follows.
\end{proof}

In the symplectic case, the above lemma is essentially \cite[Prop.2.2]{MS}.

\begin{cor} \label{nonsymb+1} Let $X$ be a symplectic $4$-manifold with $b^+=1$
containing a symplectic c-embedded torus $T$. Suppose also that
$\sw_{X\#_{T=F}E(1)}\ne0$. If $\DD_K$ is not monic, then
$X_K$ does not admit a symplectic structure.
\end{cor}
\begin{proof} Write $\DD_K(t)=a_0+\sum\limits_{j=1}^da_j(t^j+t^{-j})$
with $a_d\ne0$. Suppose that $X_K$ admits a symplectic structure with
symplectic form $\o$ and canonical class $\k$. If $[T]\cdot\o=0$, it
follows from
Lemma~\ref{lightcone} that $[T]=\lam\o$ for some $\lam\ne0$. Such a
symplectic form is
clearly nongeneric and we may perturb it so that $[T]\cdot\o\ne0$. We may
also assume
that $T$ is oriented so that  $[T]\cdot\o>0$. The adjunction inequality of
Li and Liu
\cite[Thm.E]{LL} then implies that $\k\cdot [T]\le0$.

The hypothesis that $\sw_{X\#_{T=F}E(1)}\ne0$ and Lemma~\ref{almfin} imply
that
there is a Laurent polynomial $S_X\ne0$ such that, for
$t=\exp(2[T])$, one of
\begin{gather}
\sw^-_{X,T} = S_X\cdot\sum_{n=0}^{\infty}t^{(2n+1)/2}\label{i}\\
 \sw^-_{X,T}=S_X\cdot\sum_{n=0}^{\infty}t^{-(2n+1)/2}\label{ii}
\end{gather}
holds. Let $\a$ be any class
such that the coefficient of $\exp(\a)$ in $S_X$ is nonzero. There are
finitely many such
classes; so we may list the integers $m_1<\cdots<m_r$ such that the coefficient
$S_X(\a+m_i[T])$ of $\exp(\a+m_i[T])$ in $S_X$ is nonzero. Theorem~\ref{b+1K}
implies that if the case \eqref{i} holds, then
\begin{equation}
\sw^-_{X_K}(\a+(m_1-2d+1)[T])=a_d\, S_X(\a+m_1[T])\ne0
\label{a}\end{equation}
and if the case \eqref{ii} holds then
\begin{equation}
\sw^-_{X_K}(\a+(m_r+2d-1)[T])=a_d\, S_X(\a+m_r[T])\ne0.
\label{b}\end{equation}

If $\k\cdot\o<0$ then, by a result of Liu and Ohta and Ono,
$X_K$ is the blowup of a rational or ruled surface \cite[Cor.1.4]{MS}. Every
such surface has a metric of positive scalar curvature; so it follows from the
wall-crossing formula that for each $\a$, $\sw_{X_K}^{\pm}(\a)=\pm1$ or
$0$. Equations
\eqref{a} and \eqref{b} imply in this case that $a_d=\pm1$, {\it{i.e.}} that
$\DD_K$ is monic.

Thus we may assume that $\k\cdot\o>0$. This means that $\k$ and $[T]$ both
lie in the
same component of the cone $\C_{X_K}$; so Lemma~\ref{lightcone} implies that
$\k\cdot[T]\ge0$. Since we already have the opposite inequality, we must
have $\k\cdot
[T]=0$. By Lemma~\ref{wcr},
\[ |\,\sw^-_{X_K}(\k+2m[T])\pm\sw^-_{X_K}(-\k-2m[T])\,|=1 \]
for each $m$. However $[T]\cdot\o>0$; so for $m$ large, we have
$(-\k-2m[T])\cdot\o<-\k\cdot\o$. Thus Proposition~\ref{CT} implies
that $\sw^-_{X_K}(-\k-2m[T])=0$ for
$m>>0$. This means that
$\sw^-_{X_K}(\k+2m[T])=\pm1$ for $m>>0$. But
$\sw^-_{X_K,T}=\sw^-_{X,T}\cdot\DD_K(t)$ and
$(\k+2m[T])\cdot [T]=0$; so the case \eqref{i} must hold.

Again by Proposition~\ref{CT}, we have $\sw^-_{X_K}(-\k)\ne0$; so we may write
$\k=-\a+2n[T]$ where $\sw^-_X(\a)\ne0$ and $|n|\le d$, and again we have
$m_1<\cdots<m_r$ and \eqref{a}. From \eqref{i}:
\[ \sw^-_X(\a)=\sum_{m_i \text{\,odd}<0}S_X(\a+m_i[T]). \]
But $(\a+(m_1-2d+1)[T])\cdot\o=-\k\cdot\o+(m_1+1-2(d-n))[T]\cdot\o\le-\k\cdot\o$
because $m_1<0$ (since $\sw^-_X(\a)\ne0$) and $d-n\ge0$. Thus
Proposition~\ref{CT}
implies that
$m_1=-1$ and $n=d$; so $\k=\a+2d[T]$. This means that
$\pm1=\sw^-_{X_K}(-\k)=a_d\cdot\sw^-_X(-\a)$; so
$a_d=\pm1$, {\it{i.e.}} $\DD_K$ is monic.
\end{proof}

As in the $b^+>1$ case, if $X$  contains a surface
$\Sigma_g$ of genus $g$ disjoint from $T$  with $0 \ne[\Sigma_g] \in
H_2(X;\Z)$ and with
$[\Sigma_g]^2 < 2-2g$ if $g>0$ or $[\Sigma_g]^2 < 0$ if $g=0$, then $X_K$
with the
opposite orientation  does not admit a symplectic structure. These results along
with the blowup formula of \cite{FS} imply:

\begin{cor} For any knot $K$ in $S^3$ whose Alexander polynomial $\DD_K(t)$
is not monic,
the manifolds $E(1)_K$ admit no symplectic structure with either
orientation, even after
an arbitrary number of blowups. \end{cor}

The first examples of this sort were obtained by Szabo \cite{S2}, and in
fact they are
the manifolds $E(1)_{K_k}$ where $K_k$ is the $k$-twist knot. For these
examples (with
$T=F$)
\begin{gather*}
\sw^-_{E(1)_{K_k},T}=(kt-(2k+1)+kt^{-1})\cdot\sum_{n=0}^{\infty}t^{(2n+1)/2}\\
\sw^+_{E(1)_{K_k},T}=-(kt-(2k+1)+kt^{-1})\cdot\sum_{n=0}^{\infty}t^{-(2n+1)/2}
\end{gather*}
Szabo computes the `small perturbation' invariant $\sw^0_{E(1)_{K_k}}$. As
is the case
for $E(1)$, for $m>0$, this is
$\sw^0_{E(1)_{K_k}}(m[T])=\sw^+_{E(1)_{K_k}}(m[T])$, and
for $m<0$, it is $\sw^0_{E(1)_{K_k}}(m[T])=\sw^-_{E(1)_{K_k}}(m[T])$. Hence
\[ \sw^0_{E(1)_{K_k,T}} = kt^{-1/2}-kt^{1/2},\]
as Szabo calculates.

\end{document}